%% file: docsiamart.tex
\documentclass[final,onefignum,onetabnum]{siamart220329}

\usepackage{braket,amsfonts,bm}

\usepackage{array}
\usepackage{enumerate}

\usepackage{pgfplots}
\usepackage{subcaption}

\newsiamthm{claim}{Claim}
\newsiamremark{remark}{Remark}
\newsiamremark{hypothesis}{Hypothesis}
\newsiamthm{defi}{Definition}
\newsiamthm{thm}{Theorem}
\newsiamthm{ex}{Example}
\newsiamthm{rmk}{Remark}
\newsiamthm{prop}{Proposition}
\newsiamthm{cor}{Corollary}
\newsiamthm{assum}{Assumption}
\def\citep{\cite}

\def\qed{Q.E.D.}
\crefname{hypothesis}{Hypothesis}{Hypotheses}

\usepackage[algo2e]{algorithm2e} 

\usepackage{graphicx,epstopdf}

\Crefname{ALC@unique}{Line}{Lines}

\usepackage{amsopn}

\usepackage{xspace}
\usepackage{bold-extra}
\usepackage[most]{tcolorbox}

\colorlet{texcscolor}{blue!50!black}
\colorlet{texemcolor}{red!70!black}
\colorlet{texpreamble}{red!70!black}
\colorlet{codebackground}{black!25!white!25}

\lstdefinestyle{siamlatex}{%
  style=tcblatex,
  texcsstyle=*\color{texcscolor},
  texcsstyle=[2]\color{texemcolor},
  keywordstyle=[2]\color{texemcolor},
  moretexcs={cref,Cref,maketitle,mathcal,text,headers,email,url},
}

\tcbset{%
  colframe=black!75!white!75,
  coltitle=white,
  colback=codebackground, 
  colbacklower=white, 
  fonttitle=\bfseries,
  arc=0pt,outer arc=0pt,
  top=1pt,bottom=1pt,left=1mm,right=1mm,middle=1mm,boxsep=1mm,
  leftrule=0.3mm,rightrule=0.3mm,toprule=0.3mm,bottomrule=0.3mm,
  listing options={style=siamlatex}
}

\newtcblisting[use counter=example]{example}[2][]{%
  title={Example~\thetcbcounter: #2},#1}

\newtcbinputlisting[use counter=example]{\examplefile}[3][]{%
  title={Example~\thetcbcounter: #2},listing file={#3},#1}

\DeclareTotalTCBox{\code}{ v O{} }
{ 
  fontupper=\ttfamily\color{black},
  nobeforeafter,
  tcbox raise base,
  colback=codebackground,colframe=white,
  top=0pt,bottom=0pt,left=0mm,right=0mm,
  leftrule=0pt,rightrule=0pt,toprule=0mm,bottomrule=0mm,
  boxsep=0.5mm,
  #2}{#1}

\patchcmd\newpage{\vfil}{}{}{}
\begin{tcbverbatimwrite}{tmp_\jobname_header.tex}
\title{A Projection Method for Compressible Generic Two-Fluid Model\thanks{Submitted to the editors DATE.
\funding{I acknowledge the financial support from the National Science and Technology Council of Taiwan (NSTC 112-2639-E-011-001-ASP).}}}

\author{Po-Yi Wu\thanks{Sustainable Electrochemical Energy Development (SEED) Center, National Taiwan University of Science and Technology, Taipei City 106, Taiwan (\email{r04527030@gmail.com}).}}

\headers{A Projection Method for Compressible Generic Two-Fluid Model}{Po-Yi Wu}
\end{tcbverbatimwrite}
\input{tmp_\jobname_header.tex}

\ifpdf
\hypersetup{ pdftitle={Guide to Using  SIAM'S \LaTeX\ Style} }
\fi

\begin{document}
\maketitle

\begin{tcbverbatimwrite}{tmp_\jobname_abstract.tex}
\begin{abstract}
A new projection method for a generic two-fluid model is presented in this work. Specifically, we extend the projection method, originally designed for single-phase variable density incompressible and compressible flows, to viscous compressible two-fluid flows. The key idea is that the single pressure $p$ can be uniquely determined by the products of volume fractions and densities $\phi_k \rho_k$, of the two fluids. Additionally, the stability of the method is ensured by appropriately assigning intermediate step variables at the time-discrete level and incorporating a stabilizing term at the fully discrete level. We prove the energy stability of the proposed numerical scheme, and its validity is demonstrated through three numerical tests.
\end{abstract}


\end{tcbverbatimwrite}
\input{tmp_\jobname_abstract.tex}

\section{Introduction}
Multi-fluid flows are prevalent in both natural and industrial contexts. For example, the oil/gas/water three-phase system is commonly encountered in the petroleum industry \cite{thorn2012three,jabbour1996oil,abreu2006three}; two-gas systems are employed in semiconductor manufacturing \cite{ahn1987nonlinear,laudise1998physical,merkocci2005new}; and gas/liquid systems are integral to numerous chemical engineering processes \cite{triplett1999gas,grafner2023fluidic,di2004application,whalley1987boiling}. The simplest scenario is a two-fluid system, where each fluid adheres to its own governing equations and interacts with the other fluid through a free surface or constitutive relations.

Various models for two-fluid flows have been proposed, which can be broadly categorized into two types: (i) Interface-capturing models \cite{abels2012thermodynamically,shen2015decoupled,sussman1994level,olsson2005conservative} and (ii) Averaged models \cite{bresch2010global,Ni1991,ishii2010thermo,drew1983mathematical,zuber1965}. Interface-capturing models can resolve the topology of phase distribution by either directly tracking the interface or using an order parameter to monitor the phase distribution, thereby extracting detailed physical information. However, when the phase distribution topology is complex or involves a wide range of length scales, averaged models offer a valuable compromise. While averaged models obscure interface details, they offer several advantages: (i) lower computational costs; (ii) reduced sensitivity to mesh resolution; and (iii) simpler expressions for interfacial terms.

This work performs a numerical analysis at the time-discrete, space-continuous level of the generic two-fluid model presented in \cite{bresch2010global}, which is derived from a volume-averaging procedure applied to the Navier-Stokes equations for each fluid (see \cite{ishii2010thermo} for an example). The unknowns in the model include two velocities, two densities, two volume fractions, and one pressure. Importantly, the model captures the difference in motion between the two phases, a phenomenon commonly observed in gas-liquid systems. To the best of the author's knowledge, there is no existing mathematically proven stable numerical scheme for solving this system.

To develop such a scheme, ideas from the consistent splitting method for variable-density incompressible Navier-Stokes equations \cite{guermond2000projection, li2021bound, pacheco2022consistent}, the incremental projection method for the Navier-Stokes-Allen-Cahn diffuse interface model \cite{deteix2022new}, and single-phase compressible Navier-Stokes equations \cite{gallouet2008unconditionally, grapsas2016unconditionally} are adapted. All of these methods fall within the projection method framework \cite{chorin1968numerical, guermond2006overview, bell1989second}. 

This work proposes a prediction-correction procedure within the projection method framework specifically for this two-fluid model. Unlike popular computational frameworks for single-phase or multi-phase compressible flow, such as SIMPLE, PISO, etc. \cite{Gosman1992, ISSA198640, PATANKAR19721787, ISSA1994347, Barbara2022}, the method proposed here consolidates all essential iterations into a single substep (the projection step) within the time loop, rather than iterating through all substeps within each time step. Consequently, the computational cost is now primarily driven by the iteration of the projection step.

The stability analysis follows a similar approach to that in \cite{guermond1999resultat, guermond2000projection}, but maintaining energy bounds for the potential functions presents additional challenges. To address this, auxiliary functions related to the equations of state are introduced to complete the energy estimate proof. It is demonstrated that the proposed numerical scheme is unconditionally stable at the time-discrete, space-continuous level.

At the fully discrete level (finite element discretization), applying the proposed method to the two-fluid system presents two main challenges: (i) a direct adaptation of the time-discrete scheme does not preserve mass positivity; and (ii) abnormal fluid expansion or compression may occur in the presence of discontinuities, leading to instability. These two issues are closely related, as extreme velocity divergence often occurs in regions with large mass gradients. To address these, we adopt an artificial viscosity framework \cite{kurganov2012new, nazarov2013residual} to better capture regions with large gradients. By significantly reducing abnormal expansion and compression, the artificial viscosity also helps mitigate the issue of non-positive mass.

The paper is structured as follows. Section 2 presents the governing equations, their derivation, and reformulations, establishing relationships among the unknown variables and providing a priori estimates that underpin the numerical scheme. Section 3 describes the time-discrete projection method for the two-fluid model and introduces the finite element implementation for numerical testing. Section 4 provides a stability analysis of the time-discrete scheme. Section 5 details the finite element implementation of the time-discrete space-continuous scheme proposed in Section 3. Finally, Section 6 presents three test problems solved using the finite element method, along with convergence tests to demonstrate the numerical scheme's validity.

\section{Modeling equations}
In what follows, we assume that the fluid domain $\Omega \subset \mathbb{R}^d$, where $d=2,3$, is smooth, bounded, and connected. Let $T > 0$ represent the final time. We denote the boundary of the domain by $\Gamma$. For convenience, we define $\Omega_T := \Omega \times (0,T)$ and $\Gamma_T := \Gamma \times (0,T)$. 

In this work, we develop a numerical method to compute the volume fractions $\phi_k$, densities $\rho_k$, pressure $p_k$, and velocity fields $\bm{u}_k$ for $k = g, l$, which satisfy the following averaged model for two-phase fluid flow.
\begin{subequations}\label{goveq1}
\begin{equation}
\phi_g + \phi_l = 1, \quad \text{in}~\Omega_T,
\end{equation}    
\begin{equation}
\partial_t(\phi_k \rho_k) + \nabla\cdot(\phi_k \rho_k \bm{u}_k) = 0, \quad \text{in}~\Omega_T,
\end{equation}
\begin{equation}\label{goveq1.3}
\partial_t (\phi_k \rho_k \bm{u}_k) + \nabla\cdot(\phi_k \rho_k \bm{u}_k \otimes \bm{u}_k) - \nabla\cdot(\phi_k \tau_k(\bm{u}_k)) 
+ \phi_k \nabla p_k + F_{D,k} = \phi_k \rho_k \bm{g}, \quad \text{in}~\Omega_T,
\end{equation}
where $\bm g$ is the gravitational acceleration.
To close the system, we assume
\begin{equation}\label{goveq1.4}
p_g = p_l = p, \quad \text{in}~\Omega_T.
\end{equation}
\end{subequations}
It is worth noting that if capillarity between the two fluids needs to be considered, a modified version of \eqref{goveq1.4} can be employed:
\begin{equation}\tag{2.1d'}
	p_g - p_l = F(\phi_g, \nabla\phi_g),
\end{equation}
where $F$ is a suitable function of $\phi_g$ and $\nabla\phi_g$. For the simplicity of the analysis, we consider \eqref{goveq1.4} in this work.
The stress tensor $\tau_k(\bm u_k)$ are of the form
\begin{equation}
\tau_k (\bm u_k) = 2\mu_k D(\bm u_k) + \lambda_k (\nabla\cdot\bm u_k)I,
\end{equation}
where $D(\bm u_k)$ denotes the symmetric part of the velocity gradient $\nabla\bm u_k$, $\mu_k$ are the dynamic viscosities and $\lambda_k$ are the second viscosities. The drag force $F_{D,k}$ are assumed to be of the form:
\begin{equation}\label{dragforce}
F_{D,k} = C_D |\bm u_g - \bm u_l | (\bm u_k - \bm u_{\widetilde k}),~~k=g,l,~~\widetilde{k} = l,g,
\end{equation}
where $C_D$ is assumed to be a positive function of $\phi_g$, $\phi_l$, $\rho_g$, and $\rho_l$.

\subsubsection*{Equation of state}
We assume the barotropic system such that
\begin{equation}\label{pequal}
p = \zeta_g(\rho_g) = \zeta_l(\rho_l)
\end{equation}
with
\begin{equation}\label{eos}
\zeta_g(z) = A_g z^{\gamma_g},~~~\zeta_l(z) = A_l(z^{\gamma_l} - \rho_{l,0}^{\gamma_l}) + p_{0},
\end{equation}
for some positive constants $A_g,A_l,\gamma_g,\gamma_l,\rho_{l,0},p_0$. In this work, we further assume that $\gamma_k > 1$ for $k=g,l$. We note that $\zeta_l$ follows the form known as the Tait equation (see, e.g., \cite{thompson1972compressible}), which is commonly used for the equation of state of liquids.

\subsubsection*{Boundary conditions and initial conditions}
To complete the statement of the mathematical problem, we impose that
\begin{equation}
\begin{aligned}
\rho_k|_{\Gamma_{in}} = \rho_{k,in},~~\phi_k|_{\Gamma_{in}} = \phi_{k,in},~~\bm u_k|_{\Gamma} = \bm b_k, \\
\rho_k|_{t=0} = \rho_k^0,~~\phi_k|_{t=0} = \phi_k^0,~~\bm u_k|_{t=0} = \bm u_{k,0},
\end{aligned}
\end{equation}
at all time, where $\Gamma_{in}$ is the inlet boundary such that
\begin{equation}
\Gamma_{in} := \lbrace x\in \Gamma |~ \bm b_k(x) \cdot \bm n(x) < 0 \rbrace,
\end{equation}
where $\bm n(\cdot)$ is the outer normal on $\Gamma$.

\subsection{Reformulation of the governing equations}\label{otheraspects}
Let $\alpha_k = \phi_k\rho_k$, \eqref{goveq1} becomes
\begin{subequations}
\begin{equation}\label{goveq2.1}
\frac{\alpha_g}{\rho_g} + \frac{\alpha_l}{\rho_l} = 1,~~~\text{in}~\Omega_T,
\end{equation}
\begin{equation}
\partial_t \alpha_k + \nabla\cdot(\alpha_k\bm u_k) = 0,~~~\text{in}~\Omega_T,
\end{equation}
\begin{equation}\label{goveq2.3}
\begin{aligned}
&\partial_t (\alpha_k \bm u_k) + \nabla\cdot(\alpha_k \bm u_k \otimes \bm u_k) - \nabla\cdot (\frac{\alpha_k}{\rho_k}\tau_k (\bm u_k) )\\
 +& \frac{\alpha_k}{\rho_k}\nabla p + C_D(\alpha_g,\alpha_l)|\bm u_g - \bm u_l|(\bm u_k - \bm u_{\widetilde k}) = \alpha_k \bm g,~~k=g,l,~~\widetilde{k}=l,g,~~~\text{in}~\Omega_T.
 \end{aligned}
\end{equation}
\end{subequations}
The primary reason for considering the above reformulation is to reduce the number of variables. In the following, we will show that with given $\alpha_g$ and $\alpha_l$, the quantities $\phi_g$, $\phi_l$, $\rho_g$, $\rho_l$, and 
$p$ can be uniquely determined. This justifies the expression of $C_D$ in \eqref{goveq2.3}.
\subsubsection*{Pressure balance}
By \eqref{goveq1.4}, we have
\begin{equation}\label{soundspeed}
s_g^2 d\rho_g = s_l^2 d\rho_l,~~~\text{with~}s_k = \sqrt{\frac{d\zeta_k}{d\rho_k}(\rho_k)} .
\end{equation}
We note that \eqref{soundspeed} holds only if $\zeta_k$ have nonnegative derivatives. It shall be emphasized that \eqref{eos} satisfies this condition.
By $\alpha_k = \phi_k\rho_k$, we have
\[ d\rho_g = \frac{1}{\phi_g}(d\alpha_g - \rho_g d\phi_g),~~~d\rho_l = \frac{1}{\phi_l}(d\alpha_l + \rho_l d\phi_g) .\]
Using \eqref{soundspeed}, the differential of the gaseous phase volume fraction can be expressed as
\[
d\phi_g = \frac{\phi_l s_g^2}{\phi_l \rho_g s_g^2 + \phi_g\rho_l s_l^2} d\alpha_g - \frac{\phi_g s_l^2}{\phi_l \rho_g s_g^2 + \phi_g\rho_l s_l^2}d\alpha_l .
\]
Note that we always work with the case with $\phi_k >0$. The differential of the gaseous phase density can be expressed as
\[
d\rho_g = \frac{\rho_g\rho_{l}s_{l}^2}{\alpha_{l}\rho_g^2 s_g^2 + \alpha_g\rho_{l}^2 s_{l}^2}  (\rho_{l} d\alpha_g + \rho_g d\alpha_{l}).
\]
Therefore, the differential of pressure can be expressed as
\begin{equation}\label{pandalpha}
dp = C^2(\rho_l d\alpha_g + \rho_g d\alpha_l),
\end{equation}
where 
\begin{equation}
C^2 = \frac{s_l^2 s_g^2}{\phi_g\rho_l s_l^2 + \phi_l \rho_g s_g^2} .
\end{equation}

\subsubsection*{Relations among $\phi_k$, $\rho_k$, and $\alpha_k$}
We recall \eqref{goveq2.1}, which gives us
\begin{equation}
\rho_l = \frac{\alpha_l\rho_g}{\rho_g - \alpha_g} .
\end{equation}
By the pressure equilibrium assumption \eqref{goveq1.4}, we have
\begin{equation}\label{varphi}
\varphi(\rho_g) := \zeta_g(\rho_g) - \zeta_l \left( \frac{\alpha_l\rho_g}{\rho_g - \alpha_g} \right) = 0 .
\end{equation}
Differentiating $\varphi$ with respect to $\rho_g$ yields
\[ \varphi'(\rho_g) = s_g^2 + s_l^2 \frac{\alpha_l\alpha_g}{(\rho_g - \alpha_g)^2} .\]
Therefore $\varphi$ is a non-decreasing function of $\rho_g$. For non-degenerate case $\phi_k<1$, we look for $\rho_g \in (\alpha_g,+\infty)$. Since $\varphi((\alpha_g,+\infty)) = (-\infty,+\infty)$, $\rho_g = \rho_g(\alpha_g,\alpha_l)$ is uniquely determined by solving $\varphi(\rho_g) = 0$ with given $\alpha_g$ and $\alpha_l$. Finally, $\rho_l$ and $\phi_k$, $k=g,l$ are given by
\begin{equation}\label{otherq}
\phi_g(\alpha_g,\alpha_l) = \frac{\alpha_g}{\rho_g(\alpha_g,\alpha_l)},~~\phi_l(\alpha_g,\alpha_l) = 1-\frac{\alpha_g}{\rho_g(\alpha_g,\alpha_l)}, ~~ \rho_l(\alpha_g,\alpha_l) = \frac{\alpha_l\rho_g(\alpha_g,\alpha_l)}{\rho_g(\alpha_g,\alpha_l)-\alpha_g} .
\end{equation}

\subsection{A priori estimates to \eqref{goveq2.1}-\eqref{goveq2.3}}
In this work, particular attention is paid to stability issues, and a discrete version associated with the numerical scheme is desired. Let us recall the a priori estimates associated with problem \eqref{goveq2.1}-\eqref{goveq2.3} with zero forcing terms and zero velocities on $\Gamma$ at all time\cite{bresch2010global}. The following estimates demonstrate the positivity of the mass, mass conservation, and energy stability, respectively. For $k = g,l$, we have
\begin{itemize}
\item Positivity of the mass
\begin{equation}\label{apriori1}
\alpha_k (x,t) > 0,~~~\forall (x,t)\in\Omega_T . 
\end{equation}
\item Mass conservation
\begin{equation}\label{apriori2}
\int_\Omega \alpha_k (x,t) dx = \int_\Omega \alpha_k (x,0) dx,~~~\forall t\in (0,T) .
\end{equation}
\item Energy stability
\begin{equation}\label{apriori3}
\begin{aligned}
&\sum_{k=g,l} \left[  \frac{1}{2} \frac{d}{dt} \int_\Omega \alpha_k(x,t)\bm u_k(x,t)^2 dx  + \frac{d}{dt} \int_\Omega \alpha_k(x,t) e_k (\rho_k(x,t))dx \right. \\
& \left.+ \int_\Omega \phi_k(x,t) \tau_k(\bm u_k(x,t)):\nabla\bm u_k (x,t) dx  \right] + \int_\Omega C_D  |\bm u_g(x,t) - \bm u_l(x,t)|^3 dx = 0 .
\end{aligned}
\end{equation}
\end{itemize}
In the above, $e_k(\cdot)$ is the potential energy dervied from the equation of state such that
\begin{equation}
e'_k (z) = \frac{\zeta_k(z)}{z^2} .
\end{equation}
We may choose proper constants $\rho_{k,ref}$ so that the following expression makes sense:
\begin{equation}
e_k(z) = \int_{\rho_{k,ref}}^z \frac{\zeta_k (s)}{s^2}ds .
\end{equation}

\section{A time-discrete and space-continuous numerical scheme}
The main idea behind the fractional step projection method for incompressible Navier-Stokes equations is to split the viscosity from the incompressibility, aiming to predict an intermediate velocity and then correct it with the proper pressure. Building on this concept, we propose a prediction-correction procedure that creates intermediate steps for all unknown variables and corrects them through a projection step. The corresponding steps for the time-discrete scheme are listed and explained as follows:

\subsubsection*{Prediction}
\begin{itemize}
\item Step 1: Given the previous step $\alpha_k^m$ and $\bm u_k^m$, the first step is to predict an intermediate $\alpha_k$ for time step $m+1$, say $\widetilde{\alpha}_k^{m+1}$:
\begin{equation}\label{scheme1}
\frac{\widetilde{\alpha}_k^{m+1} - \alpha_k^m}{\delta t} + \nabla\cdot(\widetilde{\alpha}_k^{m+1}\bm u_k^m) = 0 .
\end{equation}
In view of the the discussion in Section 2.2, the intermediate densities $\widetilde{\rho}_k^{m+1}$ and intermediate volume fractions $\widetilde\phi_k^{m+1}$ can be obtained from $\widetilde\alpha_k^{m+1}$. These intermediate variables are used in the velocity prediction. The prediction of $\widetilde{\alpha}_k^{m+1}$ is crucial due to the nonlinearity of the system of equations \eqref{goveq1}. Additionally, the link between another intermediate step velocities $\overline{\bm u}_k^{m+1}$ and the end-of-step velocities $\bm u_k^{m+1}$ is also established through $\widetilde{\alpha}_k^{m+1}$.

\item Step 2: Given $\widetilde{\alpha}_k^{m+1}$ from Step 1, the intermediate volume fractions $\widetilde{\phi}_k^{m+1}$ and densities $\widetilde{\rho}_k^{m+1}$ are solved by
\begin{equation}\label{scheme2}
\begin{cases}
\widetilde{\phi}_g^{m+1} + \widetilde{\phi}_l^{m+1} = 1 , \\
\widetilde{\phi}_k^{m+1}\widetilde{\rho}_k^{m+1} = \widetilde{\alpha}_k^{m+1},~~k=g,l , \\
\zeta_g(\widetilde{\rho}_g^{m+1}) = \zeta_l (\widetilde{\rho}_l^{m+1}) .
\end{cases}
\end{equation}
Indeed, the algebraic equation $\varphi(\widetilde{\rho}_g^{m+1}) = 0$ at each point $x\in \Omega$ (see Eq.\eqref{varphi}) can be solved by using a root-finding technique. Once $\widetilde{\rho}_g^{m+1}$ is determined, other quantities such as $\widetilde{\rho}_l^{m+1}$, $\widetilde\phi_g^{m+1}$, and $\widetilde{\phi}_l^{m+1}$ can be easily obtained by \eqref{otherq}.
\item Step 3 (Renormalization): Given the previous step pressure $p^m$, intermediate volume fractions $\widetilde{\phi}_k^m$, and intermediate densities $\widetilde{\rho}_k^m$, along with the intermediate volume fractions at the current step $\widetilde{\phi}_k^{m+1}$, and intermediate densities at the current step $\widetilde{\rho}_k^{m+1}$, the intermediate pressure for each phase $\widetilde{p}_k^{m+1}$ is solved by
\begin{equation}\label{scheme3}
\nabla\cdot\left( \frac{\widetilde{\phi}_k^{m+1}}{\widetilde{\rho}_k^{m+1}}\nabla\widetilde{p}_k^{m+1}  \right) = \nabla\cdot\left( \sqrt{\frac{\widetilde{\phi}_k^{m+1}\widetilde{\phi}_k^m}{\widetilde{\rho}_k^{m+1}\widetilde{\rho}_k^m}}\nabla p^m   \right) ,~~k = g,l,
\end{equation}
with the boundary constraint:
\begin{equation}\label{scheme3.1}
\frac{\widetilde{\phi}_k^{m+1}}{\widetilde{\rho}_k^{m+1}} \frac{\partial\widetilde{p}_k^{m+1}}{\partial n} = \sqrt{\frac{\widetilde{\phi}_k^{m+1}\widetilde{\phi}_k^m}{\widetilde{\rho}_k^{m+1}\widetilde{\rho}_k^m}}\frac{\partial p^m}{\partial n},~~~\text{on}~\Gamma ,~~k = g,l.
\end{equation}
Instead of rendering an intermediate step pressure from the given $\widetilde{\rho}_k$ via the equation of state, two virtual intermediate pressures are created by "rearranging" the previous step pressure $p^m$. The purpose of this is to establish a connection between the previous step pressure $p^m$ and the end-of-step pressure $p^{m+1}$. This ensures that the gradient of $p^{m+1}$ is bounded by the gradient of $p^m$, as will be discussed in Section 4.

\item Step 4. Given previous step velocities $\bm u_k^m$, and other variables obtained from Step 1 to Step 3, the intermediate velocities $\widetilde{\bm u}_k^{m+1}$ are determined by solving the following system of equations::
\begin{equation}\label{scheme4}
\begin{cases}
\begin{aligned}
\frac{\widetilde{\alpha}_g^{m+1}\widetilde{\bm u}_g^{m+1} - \alpha_g^m \bm u_g^m}{\delta t} + \nabla\cdot(\widetilde{\alpha}_g^{m+1}\bm u_g^m\otimes \widetilde{\bm u}_g^{m+1}) + \widetilde{\phi}_g^{m+1}\nabla\widetilde{p}_g^{m+1}
 - \nabla\cdot(\widetilde{\phi}_g^{m+1}\tau_g(\widetilde{\bm u}_g^{m+1})) \\
 + C_D(\widetilde{\alpha}_g^{m+1}, \widetilde \alpha_l^{m+1})|\bm u_g^m - \bm u_l^m |(\widetilde{\bm u}_g^{m+1} - \widetilde{\bm u}_{l}^{m+1}) = \widetilde{\alpha}_g^{m+1}\bm g , 
\end{aligned}\\
\begin{aligned}
\frac{\widetilde{\alpha}_l^{m+1}\widetilde{\bm u}_l^{m+1} - \alpha_l^m \bm u_l^m}{\delta t} + \nabla\cdot(\widetilde{\alpha}_l^{m+1}\bm u_l^m\otimes \widetilde{\bm u}_l^{m+1}) + \widetilde{\phi}_l^{m+1}\nabla\widetilde{p}_l^{m+1}
 - \nabla\cdot(\widetilde{\phi}_l^{m+1}\tau_l(\widetilde{\bm u}_l^{m+1})) \\
 + C_D(\widetilde{\alpha}_g^{m+1}, \widetilde \alpha_l^{m+1})|\bm u_g^m - \bm u_l^m |(\widetilde{\bm u}_l^{m+1} - \widetilde{\bm u}_{g}^{m+1}) = \widetilde{\alpha}_l^{m+1}\bm g , 
\end{aligned}
\end{cases}
\end{equation}
This step is standard within the framework of the projection method. The forcing term is included here to provide a better initial guess for the subsequent projection step.
\end{itemize}

\subsubsection*{Correction}
\begin{itemize}
\item Step 5 (Projection): Given all variables obtained from Step 1 to Step 4, another intermediate velocities $\overline{\bm u}_k^{m+1}$, and end-of-step pressure $p^{m+1}$, volume fractions $\phi_k^{m+1}$, and densities $\rho_k^{m+1}$ are solved by the following system of equations.
\begin{equation}\label{scheme5}
\begin{cases}
\displaystyle\widetilde{\alpha}_k^{m+1}\frac{\overline{\bm u}_k^{m+1} - \widetilde{\bm u}_k^{m+1}}{\delta t} + \widetilde \phi_k^{m+1}(\nabla p^{m+1} - \nabla\widetilde{p}_k^{m+1}) = 0, \\
\displaystyle \frac{\alpha_k^{m+1} - \alpha_k^m}{\delta t} + \nabla\cdot(\widetilde{\phi}_k^{m+1}\rho_k^{m+1} \overline{\bm u}_k^{m+1}) = 0 , \\
\phi_g^{m+1} + \phi_l^{m+1} = 1 , \\
\alpha_k^{m+1} = \phi_k^{m+1}\rho_k^{m+1}, \\
\rho_k^{m+1} = \zeta_k^{-1}(p^{m+1}) .
\end{cases}
\end{equation}
The above system of equations is evidently highly nonlinear. The linearization procedure and iteration method for obtaining the approximated solution will be detailed in Section 5.2. From a numerical experimentation perspective, this design enhances stability during the iteration process when dealing with nonlinearity.

\item Step 6: Finally, given $\alpha_k^{m+1}$, $\widetilde{\alpha}_k^{m+1}$ and $\overline{\bm u}_k^{m+1}$, the end-of-step velocities $\bm u_k^{m+1}$ can be obtained by
\begin{equation}\label{scheme6}
\sqrt{\alpha_k^{m+1}}\bm u_k^{m+1} = \sqrt{\widetilde{\alpha}_k^{m+1}}\overline{\bm u}_k^{m+1},~~\alpha_k^{m+1} = \phi_k^{m+1}\rho_k^{m+1} .
\end{equation}
\end{itemize}

\begin{rmk}
	As reported in \cite{guermond2000projection}, renormalization serves as a technique for the stability analysis presented in Section 4. However, in all the numerical experiments conducted in this paper, renormalization shows no observable significance compared to simply taking $\widetilde p_g^{m+1}=\widetilde  p_l^{m+1} = p^m$. 
\end{rmk}

\begin{rmk}
	If Neumann boundary conditions for the momentum balance equation \eqref{goveq2.3} are assumed on some boundary $\Gamma_N$ of $\partial \Omega$, the weak formulation for the first equation in Step 5 can be written as
\[
\int_\Omega \left( \widetilde{\alpha}_k^{m+1} \frac{\overline{\bm{u}}_k^{m+1} - \widetilde{\bm{u}}_k^{m+1}}{\delta t} + \widetilde{\phi}_k^{m+1} \nabla p^{m+1} - \widetilde{\phi}_k^{m+1} \nabla \widetilde{p}_k^{m+1} \right) \cdot \bm{v} \, dx = 0,
\]
where $\bm{v}$ is sufficiently smooth and vanishes on $\partial \Omega \setminus \Gamma_N$.
\end{rmk}

\section{Stability analysis}

We denote by $\|\cdot \|_{L^p} := \|\cdot \|_{L^p(\Omega)}$ the $L^p$ norm on $\Omega$, $\| \cdot \|_{W^{k,p}} := \|\cdot \|_{W^{k,p}(\Omega)}$ the $W^{k,p}$ norm, and $ \|\cdot\|_k = \|\cdot \|_{H^k} = \|\cdot \|_{W^{k,2}}$, $0\leq k \leq +\infty$, $1\leq p \leq +\infty$. We denote by $(\cdot,\cdot)$ the inner product associated to $L^2(\Omega)$ space such that $(u,v) = \int_\Omega u(x)v(x) dx$ for $u,v$ in $L^2(\Omega)$.

\subsubsection*{Mass conservation}
Step 5 and Step 6 guarantee the mass conservation in the following sense:
\begin{prop}
If $\overline{\bm u}_k^{m+1}$, $k=g,l$ satisfy the boundary conditions $\overline{\bm u}_k^{m+1}\cdot\bm n |_{\Gamma} = 0$, then we have
\begin{equation}
\int_\Omega \alpha_k^{m+1} dx = \int_\Omega \alpha_k^{m}  dx,~~~k=g,l .
\end{equation}
\emph{Proof.} By \eqref{scheme6} and \eqref{scheme5}, we have
\[ \frac{\alpha_k^{m+1}-\alpha_k^m}{\delta t} + \nabla\cdot(\widetilde{\phi}_k^{m+1}\rho_k^{m+1}\overline{\bm u}_k^{m+1}) = 0\]
Taking the integration on both sides, the divergence theorem together with the assumption  $\overline{\bm u}_k^{m+1}\cdot\bm n  |_{\Gamma} = 0$ lead to the conclusion. \qed
\end{prop}

\subsubsection*{Energy estimates of \eqref{scheme1}-\eqref{scheme6}}
Similar to the energy stability result given by \eqref{apriori3}, the following theorem can be asserted:
\begin{thm}\label{thm1}
For any $\delta t>0$, the solution $(\phi_k^m,\rho_k^m,\bm u_k^m, p^m)$, $m=1,2,\ldots$ of the semi-discrete scheme \eqref{scheme1}-\eqref{scheme6} with vanishing forcing term satisfies the stability estimate
\begin{equation}
\begin{aligned}
&\sum_k \left[  \frac{1}{2}\|\sqrt{\alpha_k^{m+1}}\bm u_k^{m+1} \|_0^2  + \int_\Omega \alpha_k^{m+1} e_k(\rho_k^{m+1})dx + \delta t \sum_{j=0}^m(\widetilde{\phi}_k^{j+1}\tau_k(\widetilde{\bm u}_k^{j+1}) , \nabla\widetilde{\bm u}_k^{j+1}) \right. \\ &
\left. +\frac{1}{2}\delta t^2 \left\|  \sqrt{\frac{\widetilde{\phi}_k^{m+1}}{\widetilde{\rho}_k^{m+1}}}\nabla p^{m+1}  \right\|_0^2  \right] + \delta t \sum_{j=0}^m\int_\Omega C_D (\widetilde \alpha_g^{j+1}, \widetilde \alpha_l^{j+1})|\bm u_g^{j} - \bm u_l^{j}| |\widetilde{\bm u}_g^{j+1} - \widetilde{\bm u}_l^{j+1}|^2  dx  \\
&\leq \sum_k \left(\frac{1}{2} \|\sqrt{\alpha_k^0} \bm u_k^0 \|_0^2 + \int_\Omega \alpha_k^0 e_k(\rho_k^0) dx + \frac{1}{2}\delta t^2   \left\|  \sqrt{\frac{\widetilde{\phi}_k^{0}}{\widetilde{\rho}_k^{0}}}\nabla p^{0}  \right\|_0^2  \right) .
\end{aligned}
\end{equation}
\end{thm}

\emph{Proof.} Let us begin with the estimates for the intermediate pressures $\widetilde{p}_k^{m+1}$. Taking the inner product of \eqref{scheme3} with $\widetilde{p}^{m+1}$,
\[
\left\|  \sqrt{\frac{\widetilde{\phi}_k^{m+1}}{\widetilde{\rho}_k^{m+1}}}\nabla\widetilde p_k^{m+1}  \right\|_0^2 = \left( \sqrt{\frac{\widetilde{\phi}_k^{m+1}}{\widetilde{\rho}_k^{m+1}}}\nabla\widetilde p_k^{m+1} ,  \sqrt{\frac{\widetilde{\phi}_k^{m}}{\widetilde{\rho}_k^{m+1}}} \nabla p^m  \right) \leq \left\|  \sqrt{\frac{\widetilde{\phi}_k^{m+1}}{\widetilde{\rho}_k^{m+1}}}\nabla\widetilde p_k^{m+1}    \right\|_0 \left\|  \sqrt{\frac{\widetilde{\phi}_k^{m}}{\widetilde{\rho}_k^{m}}}\nabla p^{m}   \right\|_0 .
\]
Therefore,
\begin{equation}\label{stab1}
\left\|  \sqrt{\frac{\widetilde{\phi}_k^{m+1}}{\widetilde{\rho}_k^{m+1}}}\nabla\widetilde p_k^{m+1}  \right\|_0 \leq \left\|  \sqrt{\frac{\widetilde{\phi}_k^{m}}{\widetilde{\rho}_k^{m}}}\nabla p^{m}  \right\|_0 .
\end{equation}
Summing the inner product of \eqref{scheme4} with $\widetilde{\bm u}_k^{m+1}$ and the inner product of \eqref{scheme1} with $-\frac{1}{2}|\widetilde{\bm u}_k^{m+1}|^2$, we have
\begin{equation}\label{stab2}
\begin{aligned}
&\frac{1}{\delta t}\|\sqrt{\widetilde{\alpha}_k^{m+1}}\widetilde{\bm u}_k^{m+1}\|_0^2 - \frac{1}{\delta t}(\alpha_k^m\bm u_k^m , \widetilde{\bm u}_k^{m+1}) + (\nabla\cdot(\alpha_k^m\bm u_k^m \otimes \widetilde{\bm u}_k^{m+1}), \widetilde{\bm u}_k^{m+1}) + \delta t(\widetilde{F}_{D,k}^{m+1},\widetilde{\bm u}_k^{m+1}) \\
&+(\widetilde{\phi}_k^{m+1}\nabla\widetilde{p}_k^{m+1},\widetilde{\bm u}_k^{m+1}) + (\widetilde{\phi}_k^{m+1}\tau_k(\widetilde{\bm u}_k^{m+1}),\nabla\widetilde{\bm u}_k^{m+1}) \\
&-\frac{1}{2\delta t}\|\sqrt{\widetilde{\alpha}_k^{m+1}}\widetilde{\bm u}_k^{m+1}\|_0^2 + \frac{1}{2\delta t}\|\sqrt{\alpha_k^m}\widetilde{\bm u}_k^{m+1}\|_0^2 - \frac{1}{2}(\nabla\cdot(\alpha_k^m\bm u_k^m), |\widetilde{\bm u}_k^{m+1}|^2) = 0.
\end{aligned}
\end{equation}
where $\widetilde{F}_{D,k}^{m+1} = C_D (\widetilde \alpha_g^{m+1},\widetilde \alpha_l^{m+1}) |\bm u_g^{m}-\bm u_l^{m}| (\widetilde{\bm u}_k^{m+1}-\widetilde{\bm u}_{\widetilde{k}}^{m+1})$, $k=g,l$, $\widetilde k = l,g$.
Using the Cauchy's inequality, we have 
\begin{equation}\label{stab3.1}
\frac{1}{\delta t}(\alpha_k^{m}\bm u_k^m,\widetilde{u}_k^{m+1}) \leq \frac{1}{2\delta t}\| \sqrt{\alpha_k^m}\bm u_k^m \|_0^2 + \frac{1}{2\delta t} \| \sqrt{\alpha_k^m}\widetilde{\bm u}_k^{m+1} \|_0^2.
\end{equation}
By the formula for divergence of tensor product, we have
\begin{equation}\label{stab3.2}
\nabla\cdot(\widetilde{\alpha}_k^{m+1}\bm u_k^m \otimes \widetilde{\bm u}_k^{m+1}) = \nabla\cdot(\widetilde{\alpha}_k^{m+1}\bm u_k^m)\widetilde{\bm u}_k^{m+1} + (\widetilde{\alpha}_k^{m+1}\bm u_k^m \cdot\nabla)\widetilde{\bm u}_k^{m+1}.
\end{equation}
At this point, let us recall the following relation, which is commonly used:
For sufficiently smooth $\varphi$ and $\bm v$ with $\bm v \cdot\bm n|_{\Gamma} = 0$, we have
\begin{equation}\label{lemma1}
 \int_\Omega \left(  \varphi\bm v\cdot\nabla\varphi + \frac{1}{2}\varphi^2\nabla\cdot\bm v \right) dx = 0 .
\end{equation}
In view of \eqref{stab3.2} and using an anologue of \eqref{lemma1}, we have
\begin{equation}\label{stab3.3}
\begin{aligned}
&(\nabla\cdot(\widetilde{\alpha}_k^{m+1}\bm u_k^m \otimes \widetilde{\bm u}_k^{m+1}), \widetilde{\bm u}_k^{m+1}) - \frac{1}{2}(\nabla\cdot(\widetilde{\alpha}_k^{m+1}\bm u_k^m),|\widetilde{\bm u}_k^{m+1}|^2) \\
&= \int_\Omega \nabla\cdot(\widetilde{\alpha}_k^{m+1}\bm u_k^m)|\widetilde{\bm u}_k^{m+1}|^2 dx + \int_\Omega \widetilde{\bm u}_k^{m+1} \cdot [(\widetilde{\alpha}_k^{m+1}\bm u_k^m \cdot\nabla)\widetilde{\bm u}_k^{m+1}] dx - \frac{1}{2}\int_\Omega \nabla\cdot(\widetilde{\alpha}_k^{m+1}\bm u_k^m)|\widetilde{\bm u}_k^{m+1}|^2 dx \\
&= \int_\Omega \widetilde{\bm u}_k^{m+1} \cdot [(\widetilde{\alpha}_k^{m+1}\bm u_k^m \cdot\nabla)\widetilde{\bm u}_k^{m+1}] dx + \frac{1}{2}\int_\Omega \nabla\cdot(\widetilde{\alpha}_k^{m+1}\bm u_k^m)|\widetilde{\bm u}_k^{m+1}|^2 dx = 0.
\end{aligned}
\end{equation}
Inserting \eqref{stab3.1} and \eqref{stab3.3} into \eqref{stab2}, we arrive at
\begin{equation}\label{stab3}
\begin{aligned}
\frac{1}{2}\| \sqrt{\widetilde{\alpha}_k^{m+1}}\widetilde{\bm u}_k^{m+1}\|_0^2 + \delta t(\widetilde{\phi}_k^{m+1}\nabla\widetilde{p}_k^{m+1} , \widetilde{\bm u}_k^{m+1}) + \delta t(\widetilde{\phi}_k^{m+1}\tau_k(\widetilde{\bm u}_k^{m+1}),\nabla\widetilde{\bm u}_k^{m+1} ) \\
 + \delta t(\widetilde{F}_{D,k}^{m+1},\widetilde{\bm u}_k^{m+1})  \leq \frac{1}{2} \|\sqrt{\alpha_k^m}\bm u_k^m \|_0^2 .
\end{aligned}
\end{equation}
Taking the inner product of the first equation in \eqref{scheme5} with $\delta t\overline{\bm u}_k^{m+1}$,
\begin{equation}\label{stab4}
\begin{aligned}
\frac{1}{2}\|\sqrt{\widetilde{\alpha}_k^{m+1}}\overline{\bm u}_k^{m+1}\|_0^2 + \frac{1}{2}\|\sqrt{\widetilde{\alpha}_k^{m+1}}(\overline{\bm u}_k^{m+1} - \widetilde{\bm u}_k^{m+1})\|_0^2 - \frac{1}{2}\|\sqrt{\widetilde{\alpha}_k^{m+1}}\widetilde{\bm u}_k^{m+1}\|_0^2  \\
+\delta t(\widetilde\phi_k^{m+1}\nabla p^{m+1},\overline{\bm u}_k^{m+1}) - \delta t(\widetilde{\phi}_k^{m+1}\nabla \widetilde{p}_k^{m+1},\overline{\bm u}_k^{m+1} )= 0,
\end{aligned}
\end{equation}
Taking the inner product of the first equation in \eqref{scheme5} again with $\displaystyle\delta t\frac{\nabla\widetilde{p}_k^{m+1}}{\widetilde{\rho}_k^{m+1}}$,
\begin{equation}\label{stab5}
\begin{aligned}
&(\widetilde{\phi}_k^{m+1}\nabla\widetilde{p}_k^{m+1},\overline{\bm u}_k^{m+1}) - (\widetilde{\phi}_k^{m+1}\nabla\widetilde{p}_k^{m+1}, \widetilde{\bm u}_k^{m+1})  \\
&-\frac{1}{2}\delta t\left( \left\|  \sqrt{\frac{\widetilde{\phi}_k^{m+1}}{\widetilde{\rho}_k^{m+1}}}\nabla\widetilde p_k^{m+1}  \right\|_0^2 + \left\|  \sqrt{\frac{\widetilde{\phi}_k^{m+1}}{\widetilde{\rho}_k^{m+1}}}\nabla(\widetilde p_k^{m+1} - p^{m+1}) \right\|_0^2 - \left\|  \sqrt{\frac{\widetilde{\phi}_k^{m+1}}{\widetilde{\rho}_k^{m+1}}}\nabla p^{m+1}  \right\|_0^2 \right) =0 .
\end{aligned}
\end{equation}
Combining \eqref{stab1}-\eqref{stab5}, we have
\begin{equation}\label{stabkinetic1}
\begin{aligned}
&\frac{1}{2}\|\sqrt{\widetilde{\alpha}_k^{m+1}}\overline{\bm u}_k^{m+1}\|_0^2  + \frac{1}{2}\| \sqrt{\widetilde{\alpha}_k^{m+1}}(\overline{\bm u}_k^{m+1} - \widetilde{\bm u}_k^{m+1})\|_0^2  + \delta t(\widetilde{\phi}_k^{m+1}\nabla p^{m+1},\overline{\bm u}_k^{m+1}) \\
&+\delta t(\widetilde{\phi}_k^{m+1}\tau_k(\widetilde{\bm u}_k^{m+1}),\nabla\widetilde{\bm u}_k^{m+1}) + \delta t(\widetilde{F}_{D,k}^{m+1},\widetilde{\bm u}_k^{m+1})  + \frac{1}{2}\delta t^2  \left\|  \sqrt{\frac{\widetilde{\phi}_k^{m+1}}{\widetilde{\rho}_k^{m+1}}}\nabla p^{m+1}  \right\|_0^2  \\
&\leq \frac{1}{2}\|\sqrt{\alpha_k^m}\bm u_k^m \|_0^2 + \frac{1}{2}\delta t^2 \left\|  \sqrt{\frac{\widetilde{\phi}_k^{m}}{\widetilde{\rho}_k^{m}}}\nabla p^{m}  \right\|_0^2 +  \frac{1}{2}\delta t^2 \left\|  \sqrt{\frac{\widetilde{\phi}_k^{m+1}}{\widetilde{\rho}_k^{m+1}}}\nabla(\widetilde p_k^{m+1} - p^{m+1}) \right\|_0^2 .
\end{aligned}
\end{equation}
At this stage, the issue is to estimate $\delta t(\widetilde\phi_k^{m+1}\nabla p^{m+1},\overline{\bm u}_k^{m+1})$. Let us introduce auxiliary functions $f_k(z) = ze_k(z)$. If $\gamma_k > 1$, then $f_k$ are $C^2$, strictly convex functions for $z>0$. Indeed, by definition, we have the following:
\[ f_k(z) = z\int_{\rho_{k,ref}}^z \frac{\zeta_k(s)}{s^2} ds,~~~z\in(0,+\infty) .\]
Differentiating $f$ twice with respect to $z$, we get
\[ f_k''(z) = \frac{\zeta_k(z)}{z^2} + \frac{\zeta_k'(z)z - \zeta_k(z)}{z^2} = \frac{\zeta_k'(z)}{z} .\]
Therefore, we have
\begin{equation}
\zeta_k'(z) = zf_k''(z) .
\end{equation}

We observe that
\begin{equation}\label{eq52}
\begin{aligned}
&{-\delta t}\int_\Omega \nabla\cdot(\widetilde\phi_k^{m+1}\rho_k^{m+1}\overline{u}_k^{m+1}) f_k'(\rho_k^{m+1}) dx = {\delta t}\int_\Omega \widetilde\phi_k^{m+1} \rho_k^{m+1}\overline{\bm u}_k^{m+1}\cdot\nabla\left[\frac{d}{d\rho} (\rho e_k(\rho) \right]_{\rho = \rho_k^{m+1}} dx \\
&={\delta t}\int_\Omega \widetilde\phi_k^{m+1}\overline{\bm u}_k^{m+1}\cdot\nabla \rho_k^{m+1} \left( \frac{d \zeta_k(\rho)}{d\rho} \right)_{\rho = \rho_{k}^{m+1}} dx= {\delta t}\int_\Omega \widetilde\phi_k^{m+1}\nabla p^{m+1}\cdot\overline{u}^{m+1} dx  .
\end{aligned}
\end{equation}
On the other hand, {by applying the second equation in \eqref{scheme5}, we have the following:}
\begin{equation}\label{eq53}
\begin{aligned}
&{-\delta t \nabla\cdot(\widetilde\phi_k^{m+1}\rho_k^{m+1}\overline{u}_k^{m+1}) f_k'(\rho_k^{m+1})= }~(\alpha_k^{m+1} - \alpha_k^m)f_k'(\rho_k^{m+1}) \\
& = \alpha_k^{m+1}f_k'(\rho_k^{m+1}) -\alpha_k^m f_k'(\rho_k^m) + \alpha_k^m(f_k'(\rho_k^m) - f_k'(\rho_k^{m+1})) \\
&=\alpha_k^{m+1}e_k(\rho_k^{m+1}) - \alpha_k^m e_k(\rho_k^m) + \phi_k^{m+1} \zeta_k({\rho_k^{m+1}}) - \phi_k^m \zeta_k(\rho_k^m) + \alpha_k^m(f_k'(\rho_k^m) - f_k'(\rho_k^{m+1})) .
\end{aligned}
\end{equation}
Now,
\begin{equation}\label{eq54}
\begin{aligned}
&\phi_k^{m+1} \zeta_k(\rho_k^{m+1}) - \phi_k^m \zeta_k(\rho_k^m) + \alpha_k^m(f_k'(\rho_k^m) - f_k'(\rho_k^{m+1})) \\
&=\zeta_k(\rho_k^{m+1})(\phi_k^{m+1} - \phi_k^m) + \phi_k^m(\zeta_k(\rho_k^{m+1}) - \zeta_k(\rho_k^m) + \rho_k^m(f_k'(\rho_k^m) - f_k'(\rho_k^{m+1})) .
\end{aligned}
\end{equation}
Since $\phi_k^m \geq 0$ and $f''(z)>0$, we have
\begin{equation}\label{eq55}
\phi_k^m\left[\zeta_k(\rho_k^{m+1}) - \zeta_k(\rho_k^m) + \rho_k^m(f_k'(\rho_k^m) - f_k'(\rho_k^{m+1}))\right] = \phi_k^m\left[ \int_{\rho_k^m}^{\rho_k^{m+1}} (z-\rho_k^m)f_k''(z)dz \right] \geq 0 .
\end{equation}
{By applying \eqref{eq55} to the last three terms on the right-hand side of \eqref{eq54} and the equations of state $\zeta_k(\rho_k^{m+1}) = p^{m+1}$, we obtain the following inequality by combining \eqref{eq52} and \eqref{eq53}:}
\begin{equation}
\delta t\int_\Omega \widetilde\phi_k^{m+1}\nabla p^{m+1}\cdot\overline{\bm u}_k^{m+1} dx \geq \int_\Omega [\alpha_k^{m+1} e_k(\rho_k^{m+1}) - \alpha_k^m e_k(\rho_k^m) + p^{m+1}(\phi_k^{m+1} - \phi_k^m)] dx.
\end{equation}
Taking the summation with respect to $k$, we have
\begin{equation}\label{stabpdivu}
\begin{aligned}
&\sum_{k}\delta t\int_\Omega \widetilde\phi_k^{m+1} \nabla p^{m+1}\cdot\overline{\bm u}_k^{m+1} dx \geq \sum_k \int_\Omega \left(\alpha_k^{m+1} e_k(\rho_k^{m+1}) - \alpha_k^m e_k(\rho_k^m) + p^{m+1}(\phi_k^{m+1} - \phi_k^m) \right) dx \\
&=\sum_{k} \int_\Omega [\alpha_k^{m+1} e_k(\rho_k^{m+1}) - \alpha_k^m e_k(\rho_k^m)] dx .
\end{aligned}
\end{equation}
Using \eqref{stabkinetic1}, \eqref{stabpdivu}, and \eqref{scheme6}, we have
\begin{equation}
\begin{aligned}
&\sum_k \left[  \frac{1}{2}\|\sqrt{\alpha_k^{m+1}}\bm u_k^{m+1} \|_0^2  + \int_\Omega \alpha_k^{m+1} e_k(\rho_k^{m+1}) dx + \delta t (\widetilde{\phi}_k^{m+1}\tau_k(\widetilde{\bm u}_k^{m+1}) , \nabla\widetilde{\bm u}_k^{m+1})  \right. \\
&\left. +\frac{1}{2}\delta t^2 \left\|  \sqrt{\frac{\widetilde{\phi}_k^{m+1}}{\widetilde{\rho}_k^{m+1}}}\nabla p^{m+1}  \right\|_0^2  \right]  + \delta t\int_\Omega C_D {(\widetilde{\alpha}_g^{m+1},\widetilde{\alpha}_l^{m+1})| \bm u_g^m - \bm u_l^m | |\widetilde{\bm u}_g^{m+1} - \widetilde{\bm u}_l^{m+1}|^2} dx \\
&\leq \sum_k \left(\frac{1}{2} \|\sqrt{\alpha_k^m} \bm u_k^m \|_0^2 + \int_\Omega \alpha_k^m e_k(\rho_k^m) dx + \frac{1}{2}\delta t^2   \left\|  \sqrt{\frac{\widetilde{\phi}_k^{m}}{\widetilde{\rho}_k^{m}}}\nabla p^{m}  \right\|_0^2  \right) .
\end{aligned}
\end{equation}
Summing over $m$, the proof is completed. \qed

\section{Finite element implementation}

For simplicity, we consider the problem with $\bm{u}_k = 0$ on $\Gamma$, $k=g,l$. The case with nonhomogeneous boundary conditions or outlets can be modified accordingly in a standard manner.

{Let $\Omega_h$ be a polygonal approximation of $\Omega$, $\Gamma_h$ be the boundary of $\Omega_h$, $\bm n$ is the unit outer normal on $\Gamma_h$, and $\mathcal{T}_h$ be a triangulation of $\Omega_h$ which collects disjoint triangles $K$ such that $\cup_{K\in \mathcal{T}_h}K\subset \Omega_h$ and $|\Omega_h \setminus \cup_{K\in \mathcal{T}_h}K| = 0$. We define $h_K$ to be the diameter of triangle $K$ and $h$ as a function of $x\in\Omega_h$ such that $h(x) = h_K$ if $x\in K$. We denote by $(\cdot,\cdot)$ the inner product associated with the $L^2(\Omega_h)$ space such that $(u,v) = \int_{\Omega_h} u(x)v(x) ,dx$ for $u,v$ in $L^2(\Omega_h)$.} Let us introduce finite element approximations $Y_h \subset H_1$ for densities $\rho_{k,h}$ (also for intermediate step $\widetilde{\rho}_{k,h}$), volume fractions $\phi_{k,h}$ (also for intermediate step $\widetilde{\phi}_{k,h}$), their products $\alpha_{k,h}$ (also for intermediate step $\widetilde{\alpha}_{k,h}$), and for pressure $p_h$ (also for intermediate step $\widetilde{p}_{k,h}$); $\bm{X}_{0,h}\subset \bm{H}_0^1$ for intermediate step velocities $\widetilde{\bm{u}}_{k,h}$, {$\overline{\bm{u}}_{k,h}$, and final step velocities $\bm{u}_{k,h}$.} {In particular, Lagrange finite elements shall be adopted because the space-discrete version of Step 1, Step 2, Step 5, and Step 6, in Section 3.1 involve in pointwise procedure.}

{The prediction of $\widetilde{\alpha}_{k}^{m+1}$ (Step 1) can be simply adapted from \eqref{scheme1}: Given $\alpha_{k,h}^m \in Y_h$ and $\bm u_{k,h}^m\in \bm X_{0,h}$, find $\widetilde{\alpha}_{k,h}^{m+1}$ such that for all $q_h \in Y_h$,
\begin{equation}\label{femstep1}
\left( \widetilde{\alpha}_{k,h}^{m+1} - \alpha_{k,h}^{m} , q_h \right) + \delta t(\nabla\cdot(\widetilde\alpha_{k,h}^{m+1} \bm u_{k,h}^m), q_h)  = 0.
\end{equation}}

At Step 2, the intermediate $\widetilde{\phi}_{k,h}^{m+1}$ and $\widetilde{\rho}_{k,h}^{m+1}$ can be obtained pointwisely by the given $\widetilde{\alpha}_{k,h}^{m+1}$ without generating errors due to spatial discretization since all of them belong to the same space $Y_h$. One may employ a root-finding technique (e.g.  Ridder's method\cite{ridders1979new}) to solve $\widetilde{\phi}_{k,h}^{m+1}$ and $\widetilde{\rho}_{k,h}^{m+1}$ pointwisely.

The renormalization (Step 3) for intermediate pressure $\widetilde{p}_k^{m+1}$ can be proceeded by the following problem: For $m\geq 0$, find $\widetilde{p}_k^{m+1} \in Y_h$ such that for all $w_h \in Y_h$,
\begin{equation}\label{femrenormalization}
\left( \frac{\widetilde{\phi}_{k,h}^{m+1}}{\widetilde{\rho}_{k,h}^{m+1}}\nabla\widetilde{p}_{k,h}^{m+1} , \nabla w_h   \right) = \left( \sqrt{\frac{\widetilde{\phi}_{k,h}^{m+1}\widetilde{\phi}_{k,h}^m}{\widetilde{\rho}_{k,h}^{m+1}\widetilde{\rho}_{k,h}^m}}\nabla p_{k,h}^m , \nabla w_h  \right)
\end{equation}
We note that the boundary terms are eliminated because of \eqref{scheme3.1}.

The weak formulation for Step 4 (advection-diffusion step) reads: For $m\geq 0$, find $\widetilde{\bm u}_{k,h}^{m+1} \in \bm X_{0,h}$ such that for all $\bm v_h \in \bm X_{0,h}$,
\begin{equation}\label{femstep4}
\begin{aligned}
&\left(  \frac{\widetilde{\alpha}_{k,h}^{m+1}\widetilde{\bm u}_{k,h}^{m+1} - \alpha_{k,h}^m\bm u_{k,h}^m}{\delta t} , \bm v_h  \right) + \left( \nabla\cdot(\widetilde{\alpha}_{k,h}^{m+1} \bm u_{k,h}^m \otimes \widetilde{\bm u}_{k,h}^{m+1}) , \bm v_h  \right)\\
& - \left(\widetilde{p}_{k,h}^{m+1},\nabla\cdot(\widetilde{\phi}_{k,h}^{m+1} \bm v_h)\right) 
+(\widetilde{\phi}_{k,h}^{m+1}\tau_k (\widetilde{\bm u}_{k,h}^{m+1}),\nabla\bm v_h)\\
&+\left( {C_D(\widetilde{\alpha}_{g,h}^{m+1},\widetilde{\alpha}_{l,h}^{m+1})}| \bm u_{g,h}^{m} - \bm u_{l,h}^{m}|(\widetilde{\bm u}_{k,h}^{m+1} - \widetilde{\bm u}_{\widetilde k,h}^{m+1}), \bm v_h  \right)  = (\widetilde{\alpha}_{k,h}^{m+1}\bm g,\bm v_h),~k=g,l,~\widetilde{k}=l,g.
\end{aligned}
\end{equation}
The above equations for $g$ and $l$ form a system and are solved simultaneously.

For Step 5, we first consider its nonlinear fully discrete scheme adapted from the time-discrete version: For $m\geq 0$, find $\alpha_{k,h}^{m+1} \in Y_h$, $\overline{\bm u}_{k,h}^{m+1} \in \bm X_{0,h}$ such that the following two equations hold simultaneously for all $q_h \in Y_h$ and for all $\bm v_h \in \bm X_{0,h}$:

\begin{equation}\label{step5.1}
	(\alpha_{k,h}^{m+1}-\alpha_{k,h}^m, q_h) + \delta t(\nabla\cdot(\widetilde\phi_{k,h}^{m+1}\rho_{k,h}^{m+1}\overline{\bm u}_{k,h}^{m+1}),q_h) + \delta t(\mathcal P_{k,h}^{m+1} \nabla\alpha_{k,h}^{m+1},\nabla q_h  ) = 0,~~~k=g,l,
\end{equation}
\begin{equation}\label{step5.2}
\begin{aligned}
	&\left(\widetilde{\alpha}_{k,h}^{m+1}(\overline{\bm u}_{k,h}^{m+1} - \widetilde{\bm u}_{k,h}^{m+1}),\bm v_h  \right) + \delta t\left( \widetilde\phi_{k,h}^{m+1}\nabla p(\alpha_{g,h}^{m+1},\alpha_{l,h}^{m+1}) -\widetilde{\phi}_{k,h}^{m+1}\nabla \widetilde p_{k,h}^{m+1}, \bm v_h   \right)\\
& + \delta t\left( \eta_{k,h}^{m+1} \nabla \cdot \overline{\bm u}_{k,h}^{m+1}, \nabla \cdot \bm v_h \right)= 0,~~~k=g,l,
\end{aligned}
\end{equation}
where $\phi_{k,h}^{m+1}$ and $p$ is a function of $\alpha_{g,h}^{m+1}$ and $\alpha_{l,h}^{m+1}$ as described in Section \ref{otheraspects}. The third term on the left-hand side of \eqref{step5.1} represents the artificial diffusion, introduced to enhance stability at the fully discrete level. This term is designed to prevent the degeneration of $\alpha_k$ caused by spurious oscillation of $\nabla \cdot \bm{u}_k$.
{For \eqref{step5.2}, $\eta_{k,h}^{m+1}$ is a stabilizer function to be determined. The last term on the left-hand side of \eqref{step5.2} is referred to as artificial bulk viscosity \cite{mani2009suitability,cook2005hyperviscosity}. Its purpose is to smooth out shocks without directly dissipating vortical motions. This term is expected to remain small away from shocks and regions with steep volume fraction gradients. We observe that instability arises in regions with a large $\nabla\cdot\bm{u}_k$, which corresponds to a spurious response to the propagation of density waves. To resolve this, we design $\mathcal P_{k,h}^{m+1}$ and $\eta_{k,h}^{m+1}$ by 
\begin{equation}\label{etak}
 \mathcal P_{k,h}^{m+1}= C_{\alpha}h^2 |\nabla\cdot\widetilde{\bm u}_{k,h}^{m+1}|,~~~\eta_{k,h}^{m+1} = C_{\eta}h^2 \widetilde \alpha_{k,h}^{m+1}| \nabla\cdot \widetilde{\bm u}_{k,h}^{m+1} |,
\end{equation}
where $C_\alpha$ and $C_\eta$ are dimensionless constants. For numerical experiments, we may choose $0.1\leq C_\alpha,~C_\eta \leq 1$.

If the time step \(\delta t\) is not sufficiently small, a straightforward Picard iteration to approximate the coupled equations \eqref{step5.1} and \eqref{step5.2} may fail to converge properly. To address this issue, we introduce a positive integer \(N\) and define a substep \(\delta t_N = \delta t / N\). The coupled equations are then solved iteratively over \(N\) substeps to approximate \(\alpha_{k,h}^{m+1}\) and \(\overline{\bm{u}}_{k,h}^{m+1}\), as shown below:

\begin{equation}\label{step5.1p}\tag{5.4'}
\begin{aligned}
	&(\alpha_{k,h}^{m+1,n+1}-\alpha_{k,h}^{m+1,n}, q_h) + \delta t_N(\nabla\cdot(\widetilde\phi_{k,h}^{m+1}\rho_{k,h}^{m+1,n+1}\overline{\bm u}_{k,h}^{m+1,n+1}),q_h)\\
 &+ \delta t_N(\mathcal P_{k,h}^{m+1} \nabla\alpha_{k,h}^{m+1,n+1},\nabla q_h  ) = 0,~~~k=g,l,~~n=0,1,\ldots,N-1.
\end{aligned}
\end{equation}
\begin{equation}\label{step5.2p}\tag{5.5'}
\begin{aligned}
	&\left(\widetilde{\alpha}_{k,h}^{m+1}(\overline{\bm u}_{k,h}^{m+1,n+1} - \overline{\bm u}_{k,h}^{m+1,n}),\bm v_h  \right) + \delta t_N\left( \widetilde\phi_{k,h}^{m+1}\nabla p(\alpha_{g,h}^{m+1,n+1},\alpha_{l,h}^{m+1,n+1}) -\widetilde{\phi}_{k,h}^{m+1}\nabla \widetilde p_{k,h}^{m+1}, \bm v_h   \right)\\
& + \delta t_N\left( \eta_{k,h}^{m+1} \nabla \cdot \overline{\bm u}_{k,h}^{m+1,n+1}, \nabla \cdot \bm v_h \right)= 0,~~~k=g,l,~~n=0,1,\ldots,N-1,
\end{aligned}
\end{equation}
where $\alpha_{k,h}^{m+1,0} = \alpha_{k,h}^m$ and $\overline{\bm u}_{k,h}^{m+1,0} = \widetilde{\bm u}_{k,h}^{m+1}$.

The solutions to the coupled equations \eqref{step5.1p} and \eqref{step5.2p} are approximated using a Picard iteration. Let us consider the quasi-linearized form of \eqref{step5.1p} and \eqref{step5.2p}:
\begin{enumerate}[(i)]
\item Given \( \overline{\bm u}_{k,h}^{m+1,n+1,*} \in \bm X_{0,h} \), find \( \alpha_{k,h}^{m+1,n+1} \in Y_h \) such that the following holds for all \( q_h \in Y_h \):
\begin{equation}\label{step5.1lin}
\begin{aligned}
	&(\alpha_{k,h}^{m+1,n+1} - \alpha_{k,h}^{m+1,n}, q_h) + \delta t \left( \nabla \cdot \left( \widetilde\phi_{k,h}^{m+1} \rho_k(\alpha_{g,h}^{m+1,n+1},\alpha_{l,h}^{m+1,n+1}) \overline{\bm u}_{k,h}^{m+1,n+1,*} \right), q_h \right) \\
& +\delta t_N( \mathcal P_{k,h}^{m+1} \nabla\alpha_{k,h}^{m+1,n+1},\nabla q_h  ) = 0, \quad k = g, l,
\end{aligned}
\end{equation}

\item Given \( \alpha_{g,h}^{m+1,n+1,*}, \alpha_{l,h}^{m+1,n+1,*} \in Y_h \), find \( \overline{\bm u}_{k,h}^{m+1,n+1} \in \bm X_{0,h} \) such that the following holds for all \( \bm v_h \in \bm X_{0,h} \):
\begin{equation}\label{step5.2lin}
\begin{aligned}
	&\left( \widetilde{\alpha}_{k,h}^{m+1}(\overline{\bm u}_{k,h}^{m+1,n+1} - \overline{\bm u}_{k,h}^{m+1,n}), \bm v_h \right) + \delta t_N \left( \widetilde\phi_{k,h}^{m+1} \nabla p(\alpha_{g,h}^{m+1,n+1,*}, \alpha_{l,h}^{m+1,n+1,*}) - \widetilde{\phi}_{k,h}^{m+1,n+1} \nabla \widetilde{p}_{k,h}^{m+1}, \bm v_h \right) \\
	& + \delta t_N\left( \eta_{k,h}^{m+1} \nabla \cdot \overline{\bm u}_{k,h}^{m+1}, \nabla \cdot \bm v_h \right) = 0, \quad k = g, l.
\end{aligned}
\end{equation}
\end{enumerate}
The quasi-linearization in (i) can be easily handled using a Picard iteration.
With the above linearization, the iteration procedure of Step 5 is described by Algorithm \ref{algo1}.
\begin{algorithm}
	\SetKwInOut{Input}{initial guess} 
	\SetKwInOut{Output}{output} 
	\Input{$\overline{\bm u}_{k,h}^{m+1,0,*} \gets \widetilde{\bm u}_{k,h}^{m+1}$, $k=g,l$}
	\Output{$\alpha_{k,h}^{m+1}$, $\overline{\bm u}_{k,h}^{m+1}$, $k=g,l$}
	\For{$0\leq n \leq N-1$}{
	\For{$j$}
	{
		solve \eqref{step5.1lin} to get $\alpha_{k,h}^{m+1,n+1}$;\\
		compute the $L^2$ error: $\epsilon_{\alpha} = \left(  \sum_{k=g,l} \|\alpha_{k,h}^{m+1,n+1} - \alpha_{k,h}^{m+1,n+1,*}\|_0^2 \right)^{1/2}$;\\
		$\alpha_{k,h}^{m+1,n+1,*} \gets  \alpha_{k,h}^{m+1,n+1}$;\\
		solve \eqref{step5.2lin} to get $\overline{\bm u}_{k,h}^{m+1,n+1}$;\\
		compute the $L^2$ error: $\epsilon_{u} = \left( \sum_{k=g,l}\| \overline{\bm u}_{k,h}^{m+1,n+1} - \overline{\bm u}_{k,h}^{m+1,n+1,*}\|_0^2 \right)^{1/2}$;\\
		\If{$(\epsilon_\alpha + \epsilon_u)^{1/2} < tolerance $}{ break; }
		$\overline{\bm u}_{k,h}^{m+1,n+1,*} \gets \overline{\bm u}_{k,h}^{m+1,n+1}$;
	}
	}
	\caption{Picard iteration for the approximation solutions of Step 5.}\label{algo1}
\end{algorithm}

Step 6 is rather straightforward; $\bm u_{k,h}^{m+1}$ is obtained pointwisely using the following equation.
\begin{equation}\label{femstep6}
	\bm u_{k,h}^{m+1} = \sqrt{\widetilde{\alpha}_{k,h}^{m+1}/\alpha_{k,h}^{m+1} }\overline{\bm u}_{k,h}^{m+1}.
\end{equation}
Now we can summarize the procedure for each time step in Algorithm \ref{algo2}.
\begin{algorithm}
	\For{time step $m$}
	{
		Step 1: solve \eqref{femstep1};\\
		Step 2: solve \eqref{scheme2} pointwisely;\\
		Step 3: solve \eqref{femrenormalization};\\
		Step 4: solve \eqref{femstep4};\\
		Step 5: run Algorithm \ref{algo1};\\
		Step 6: solve \eqref{femstep6};
	}
	\caption{Algorithm for the finite element implementation.}\label{algo2}
\end{algorithm}

\section{Numerical test}\label{numtest}
Let $L$ be the reference length scale, $U$ be the reference velocity, and $\rho_0$ be the reference density. Common parameters utilized in all numerical experiments presented in this section are listed in Table \ref{tbl:1}. The parameters and variables in use are scaled by the following:
\[
\begin{aligned}
&\mu_g \rightarrow \frac{\mu_g}{\rho_0 UL},~~\mu_l \rightarrow \frac{\mu_l}{\rho_0 UL},~~ \lambda_g \rightarrow \frac{\lambda_g}{\rho_0 UL},~~\lambda_l \rightarrow \frac{\lambda_l}{\rho_0 UL}   ,~~\rho_{l,0} \rightarrow \frac{\rho_{l,0}}{\rho_0},~~A_g \rightarrow \frac{A_g \rho_0^{\gamma_g}}{\rho_0 U^2},~~A_l \rightarrow \frac{A_l \rho_0^{\gamma_l}}{\rho_0 U^2},~~p_0 \rightarrow \frac{p_0}{\rho_0 U^2},  \\
&\alpha_g \rightarrow \frac{\alpha_g}{\rho_0},~~\alpha_l \rightarrow \frac{\alpha_l}{\rho_0},~~\bm u_g \rightarrow \frac{\bm u_g}{U},~~\bm u_l \rightarrow \frac{\bm u_l}{U},~~p \rightarrow \frac{p}{\rho_0 U^2} .
\end{aligned}
\]

\begin{table}
\centering
\begin{tabular}{c  c | c  c}
\hline
Parameters   & Values & Parameters & Values \\
\hline
$\rho_{g0}$ & $1.161~(kg/m^3)$ & $\rho_{l0}$ & $995.65~~(kg/m^3)$ \\
$\mu_g$ & $1.86\times 10^{-5}~(kg/m\cdot s)$ & $\mu_l$ & $8.88\times 10^{-4}~(kg/m\cdot s)$ \\
$\lambda_g$   & $5.3847\times 10^{-4}~(kg/m\cdot s)$   &  $\lambda_l$   &   $2.47\times 10^{-3}~(kg/m\cdot s)$  \\
$A_g$ & $8.22151\times 10^4~ (m^{3.2}/kg^{0.4}\cdot s^2)$ & $A_l$ & $6~ (m^{12.2}/kg^{3.4}\cdot s^2)$ \\
$\gamma_g $ & $1.4$ & $\gamma_l$ & $4.4$ \\
$p_0$& $ 1.01325\times 10^{5}~(kg/m\cdot s^2)$ &  & \\
\hline

\end{tabular}
\caption{Parameters utilized in all numerical experiments presented in Section 5.}
\label{tbl:1}
\end{table}
The computer program for the implementation is written using Freefem++ \cite{hecht2012new}.

\subsection{Channel flow with uniform inlet distribution of fluids}
The first test is designed to demonstrate the convergence of the numerical scheme. A $P_2$ element is used for approximating the velocities, while a $P_1$ element is used for all other variables. We consider a rectangular channel with dimensions $5 \text{ m} \times 1 \text{ m}$, where the flow enters from the left (inlet), exits from the right (outlet), and has wall boundary conditions at the top and bottom. The boundary conditions are illustrated in Figure \ref{test1domain}. The initial conditions are set as $\phi_g(0) = 0.1$, $\phi_l(0) = 0.9$, and $\bm u_l (0) = \bm u_g(0) = 0$, while $\rho_g(0)$, $\rho_l(0)$, and $p(0)$ are chosen to satisfy the hydrostatic condition. The inlet pressure is given by $5(\phi_g(0)\rho_{g0} + \phi_l(0)\rho_{l0}) + p_0$, and the outlet pressure is set to $p_0$. For the drag force, we set 
\begin{equation}\label{cdmodel}
 C_D = \displaystyle\frac{\widetilde{C}_D}{L_r} \displaystyle\frac{\alpha_g \alpha_l}{\alpha_g + \alpha_l},
\end{equation}
where \( L_r \) is the characteristic length scale of the dispersed phase in the region dominated by the other phase, and \( \widetilde{C}_D \) is the dimensionless drag coefficient. Specifically, \( L_r \) can be taken to be on the same order as the specific surface area of the dispersed phase. We note that \( L_r \) should be scaled as \( L_r \rightarrow L_r / L \). In this case, we take \( \widetilde{C}_D = 0.1 \) and \( L_r = 1 \times 10^{-3} (m) \) to simulate the bubbly flow regime with a representative radius on the order of millimeters. To ensure stability, the coefficient $C_\alpha$ and $C_\eta$ in \eqref{etak} are set to 0.5 and 1, respectively.

Despite the seeming simplicity of this test due to the uniformity of initial volume fractions, the density difference between the two phases poses computational challenges. In fact, without the additional stabilizing term presented in \eqref{step5.2}, we observe abnormal expansion and compression of fluids at the inlet and outlet, leading to simulation instability.

For the convergence test, all computations are performed using a $100 \times 20$ uniform mesh. The solution obtained with a very small time step is considered as the reference solution. The convergence results depicted in Figure \ref{test1conv} indicate a first-order convergence of the numerical scheme.
\begin{figure}
\centering
\includegraphics[width=0.9\textwidth]{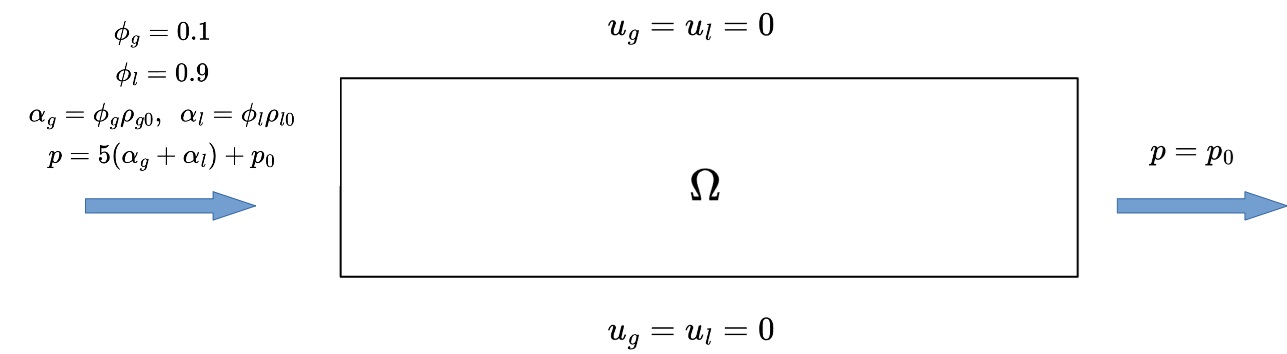}
\caption{The physical domain and boundary conditions considered in the first test.}
\label{test1domain}
\end{figure}
\begin{figure}
\includegraphics[width=\textwidth]{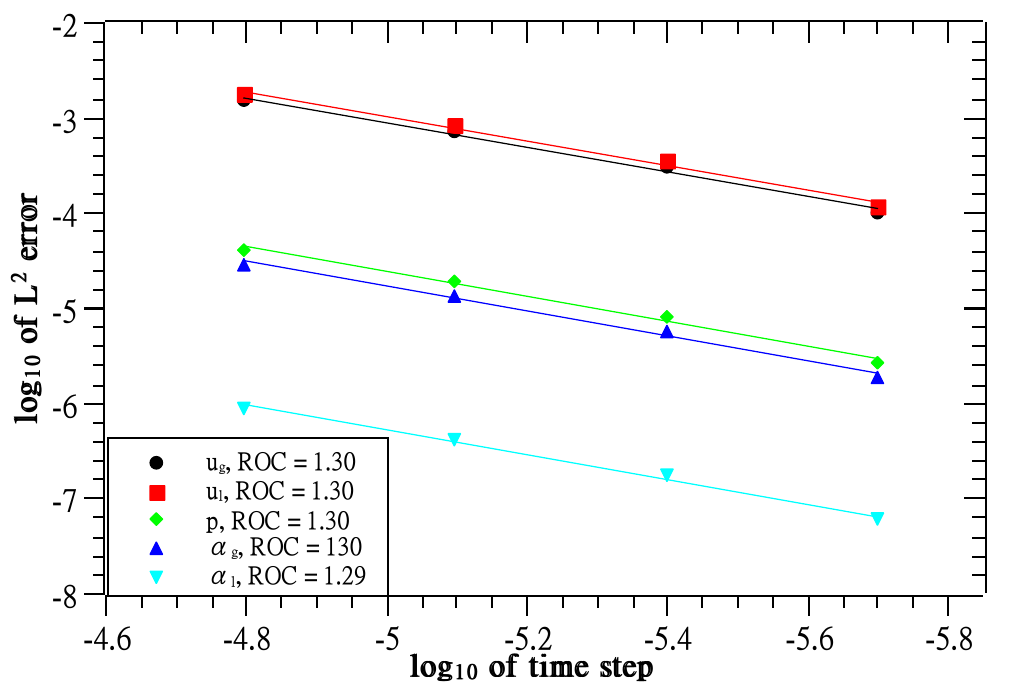}
\caption{Plot of convergence of the first test.}
\label{test1conv}
\end{figure}

\subsection{Dam break}
Second, we simulate the well-known dam break problem of Martin and Moyce \cite{Martin1952}. A $P_2$ element is used for approximating the velocities, while a $P_1$ element is used for all other variables. The physical domain is a rectangle measuring $0.5 \text{ m} \times 0.15 \text{ m}$. The test consists of a simple configuration, as shown in Figure \ref{test2domain}. Initially, a water column, $a = 0.06 \text{ m}$ wide and $\eta^2 a = 0.12 \text{ m}$ high, is at rest, and the pressure follows a hydrostatic profile. {Here, $\eta$ is a dimensionless constant for scaling the length along the $z$-axis, as depicted in Figure \ref{test2domain}.} To avoid degenerate cases where the mass of liquid or gas approaches zero, we introduce a perturbation by uniformly setting $\phi_g = 0.01$ and $\phi_l = 0.01$ in the water column and the surrounding environment, respectively. Here, a hyperbolic tangent function is employed to smooth the interface between the water column and the surrounding environment.

In this case, the drag force term is given by \eqref{cdmodel} with $\widetilde{C}_D = 1$. Since we aim to simulate the scenario where the water column region and its surrounding ambience are single-phase regions, we take a very small characteristic length scale, \( L_r = 1 \times 10^{-12} \). 
To accurately capture the sharp interface, stabilization introduced in Eqs. \eqref{step5.1} and \eqref{step5.2} should be incorporated into the fully discrete scheme. Here, we set \(C_\alpha = 0.5, C_\eta = 1 \).
All boundaries are solid walls with slip boundary conditions. Under the influence of gravity, $g = 9.8 \text{ m/s}^2$, the column collapses.

The computation is performed using a $50 \times 15$ uniform mesh. The experimental results from \cite{Martin1952} for the front position $x/a = F_1(\eta^2, t\sqrt{2g/a})$ and the height of the column $z/(\eta^2 a) = F_2(\eta^2, t\sqrt{g/a})$ are used for comparison. Note that in the experimental results, the non-dimensional times are different for the front position and the water height. We have used the same non-dimensional units in presenting the numerical results to allow a direct comparison with the original publication of \cite{Martin1952}. In Figure \ref{fig:test2}, the numerical results are compared with the experimental data for the case $\eta^2 = 2$, where $\phi_l = 0.5$ represents the interface. The numerical scheme shows good agreement with the experimental results. Figure \ref{fig:test2iso} shows the isovalues of the volume fraction at the different physical time $t=0,~ 0.066,~0.109,~0.164,~0.222,~0.281(s)$.

\begin{figure}
\centering
\includegraphics[width=0.9\textwidth]{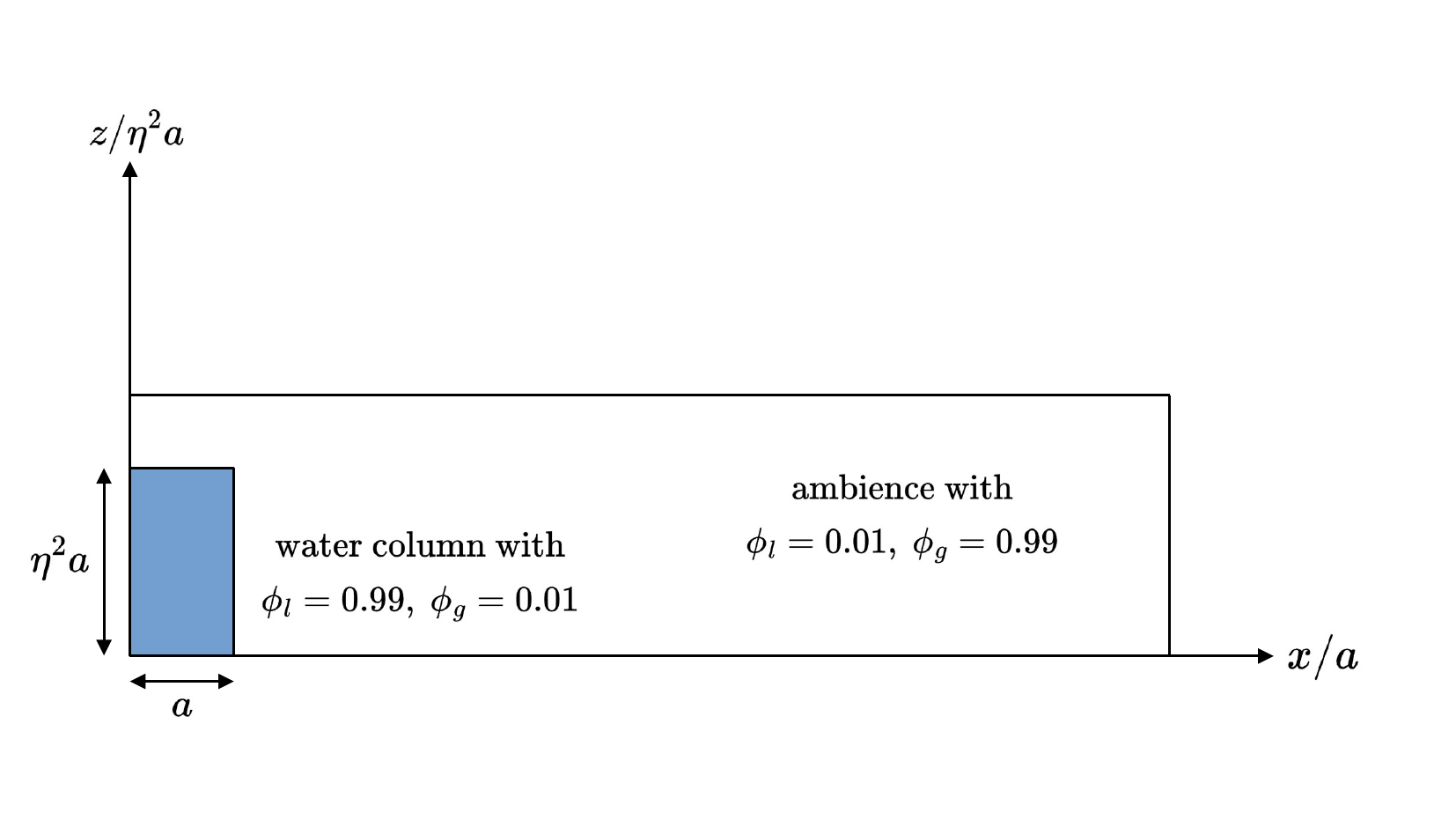}
\caption{Initial configuration of the dam break problem.}
\label{test2domain}
\end{figure}

\begin{figure}
\centering
	\begin{subfigure}{0.6\textwidth}
		\centering
		\includegraphics[width=1.0\linewidth]{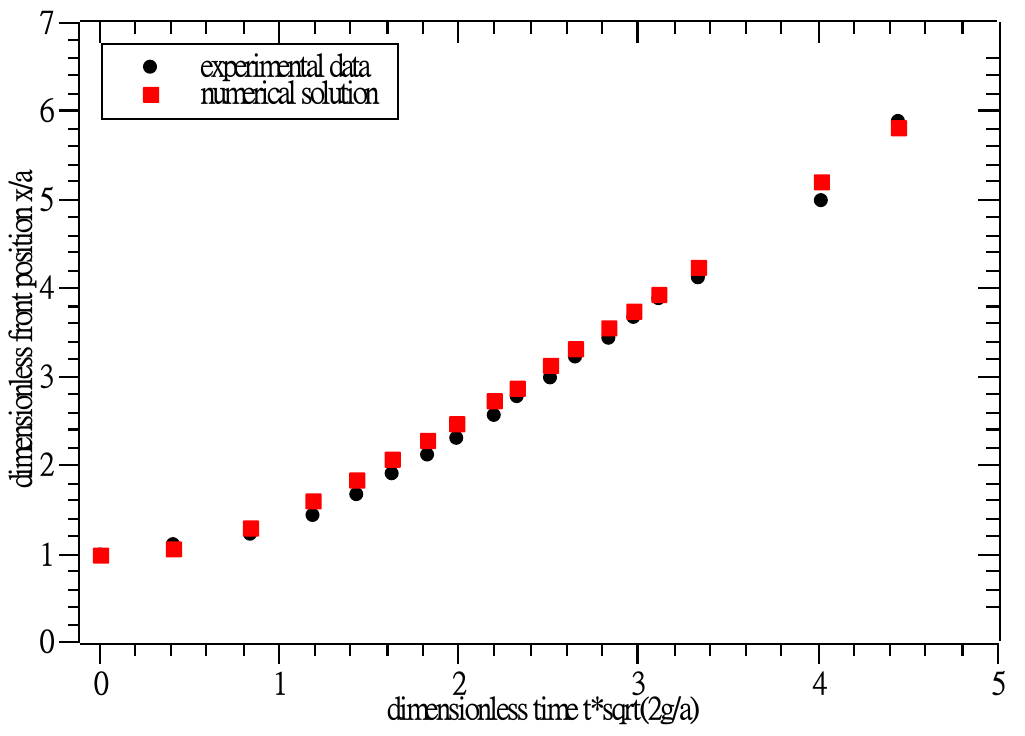}
	\end{subfigure}\hspace{\fill}%
	\begin{subfigure}{0.6\textwidth}
		\centering
		\includegraphics[width=1.0\linewidth]{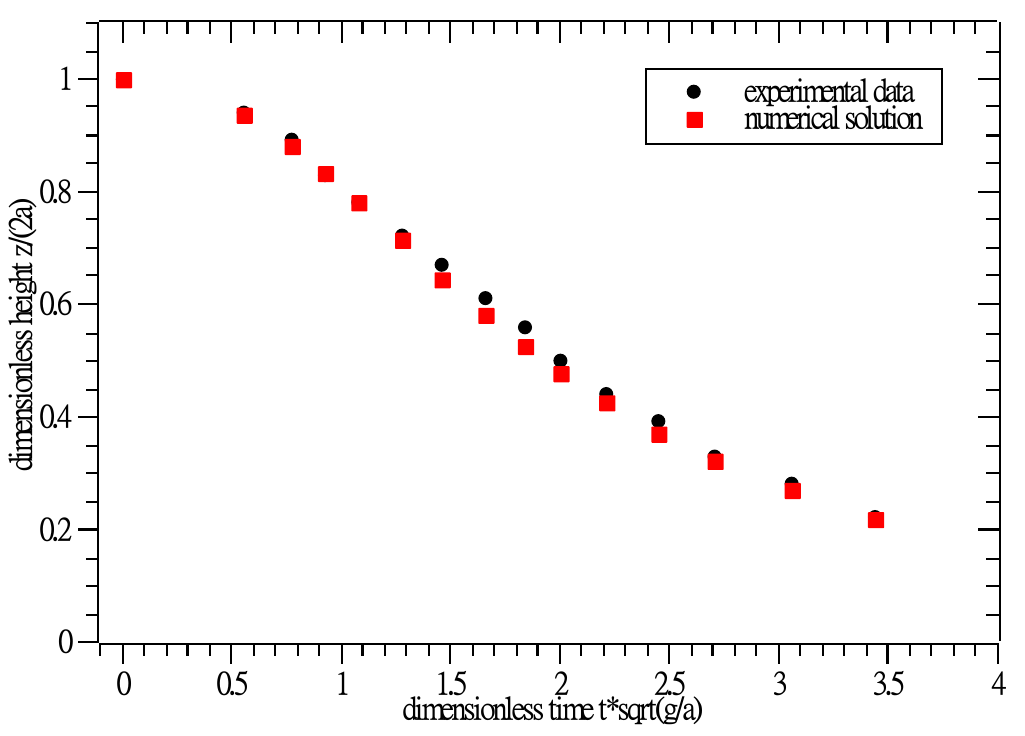}
	\end{subfigure}%
	\caption{Comparison between numerical solution and experimental results for the dam break problem. Front position (top) and height of the column (bottom).}
	\label{fig:test2}
\end{figure}

\begin{figure}
\centering
	\begin{subfigure}{0.5\textwidth}
		\centering
		\includegraphics[width=1.0\linewidth]{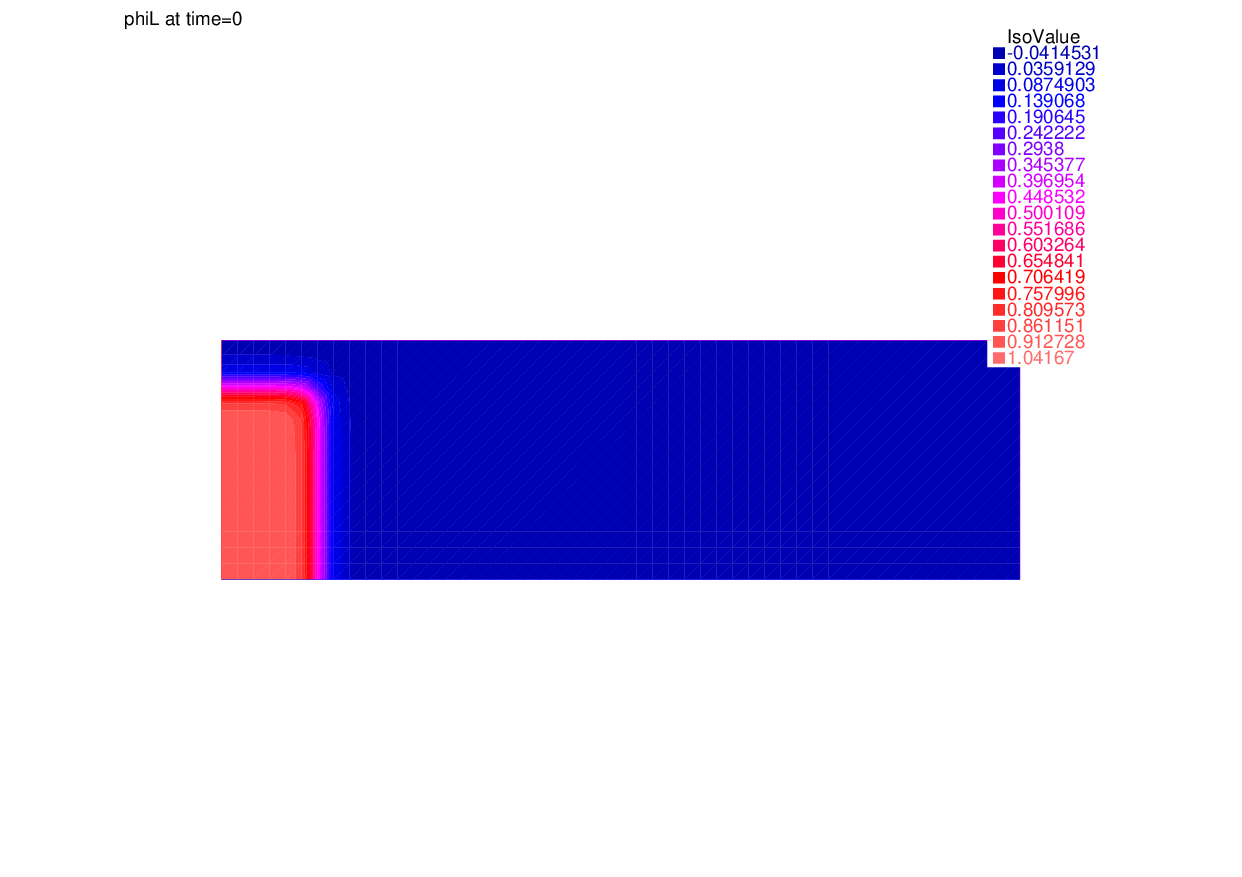}
		\caption{t=0(s)}
	\end{subfigure}\hspace{\fill}%
	\begin{subfigure}{0.5\textwidth}
		\centering
		\includegraphics[width=1.0\linewidth]{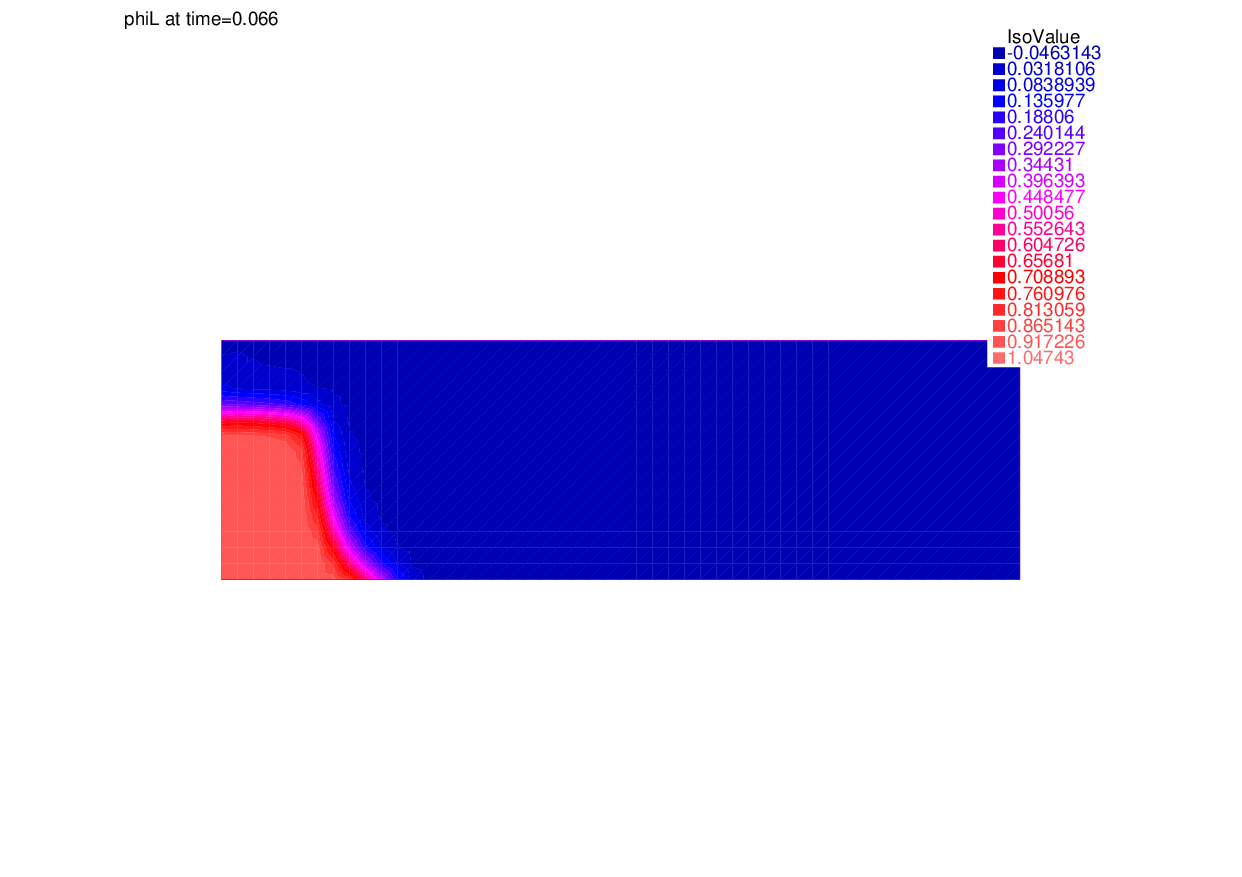}
		\caption{t=0.066(s)}
	\end{subfigure}
	\begin{subfigure}{0.5\textwidth}
		\centering
		\includegraphics[width=1.0\linewidth]{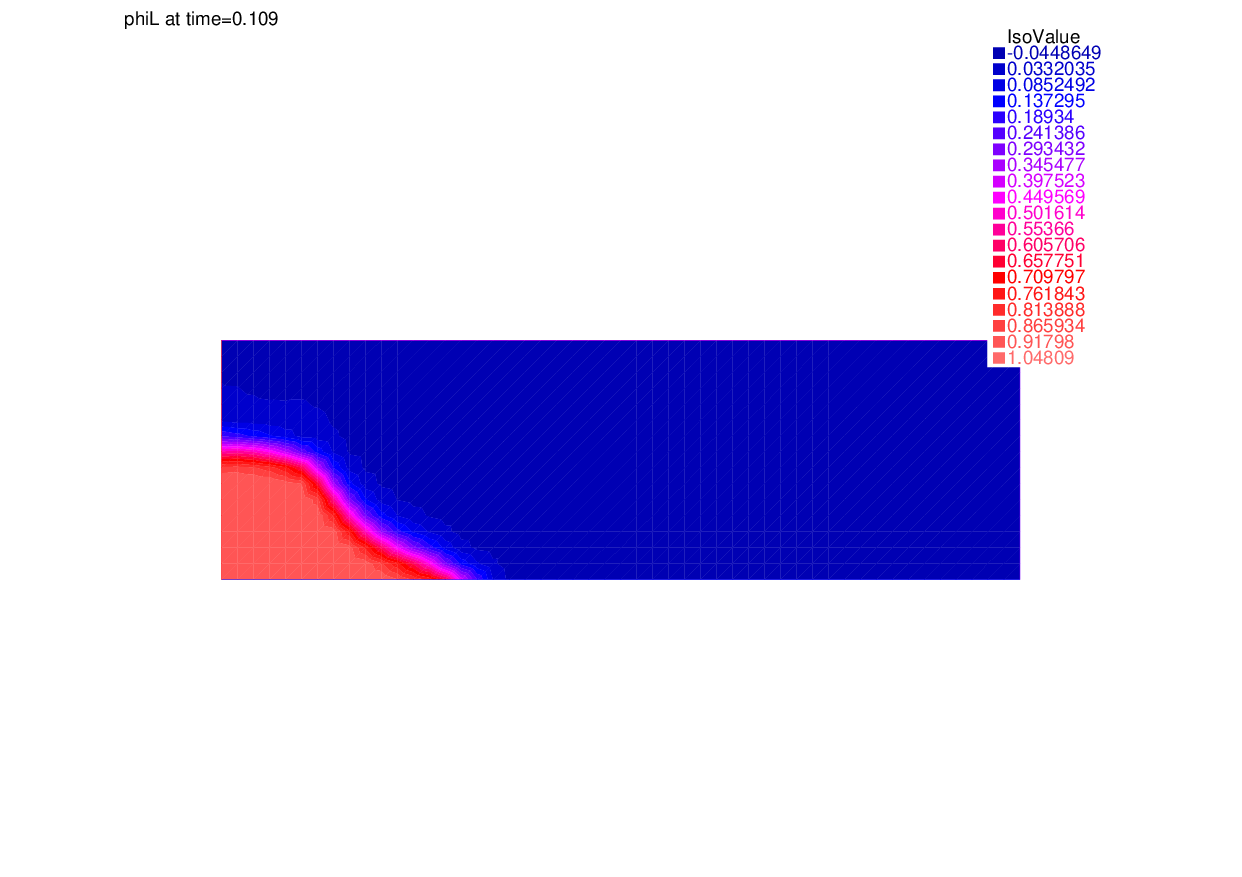}
		\caption{t=0.109(s)}
	\end{subfigure}\hspace{\fill}%
	\begin{subfigure}{0.5\textwidth}
		\centering
		\includegraphics[width=1.0\linewidth]{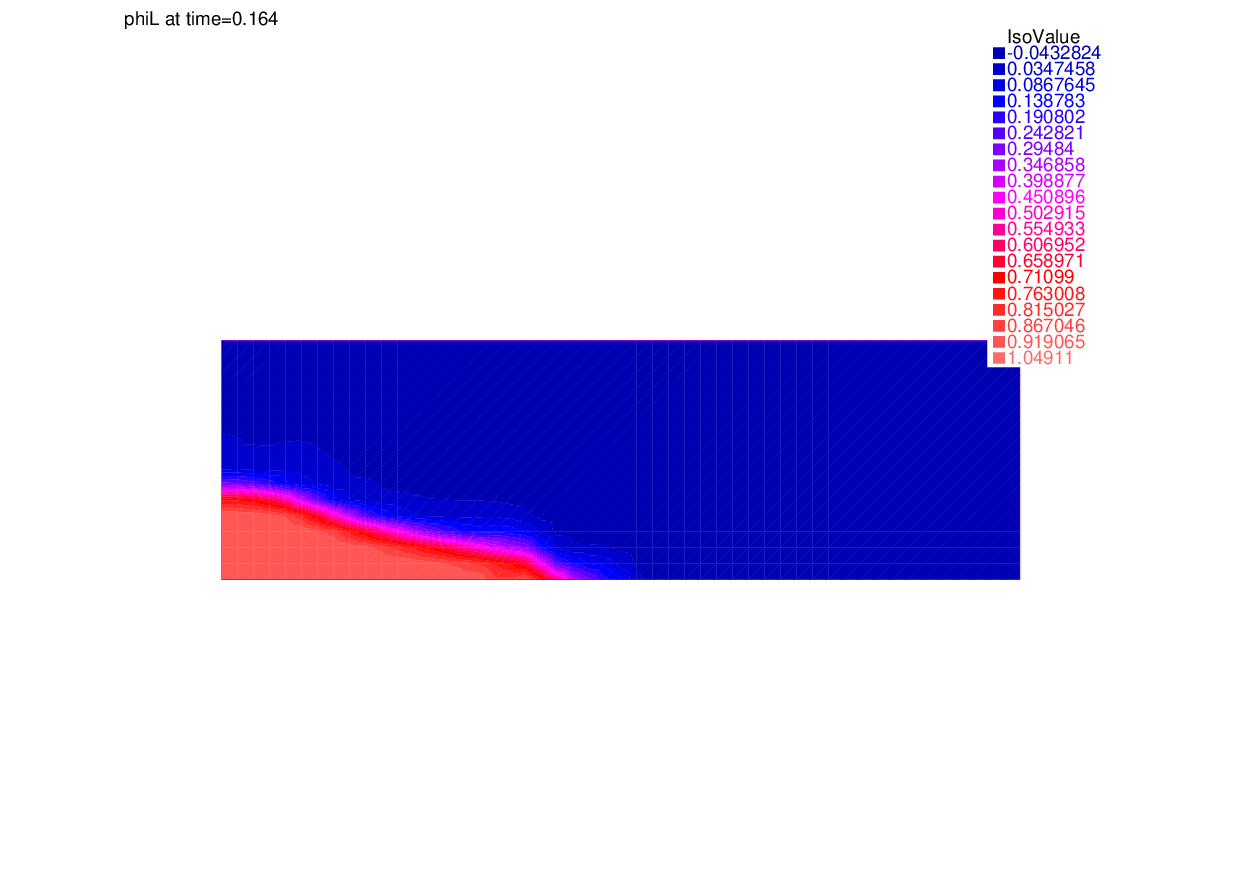}
		\caption{t=0.164(s)}
	\end{subfigure}
	\begin{subfigure}{0.5\textwidth}
		\centering
		\includegraphics[width=1.0\linewidth]{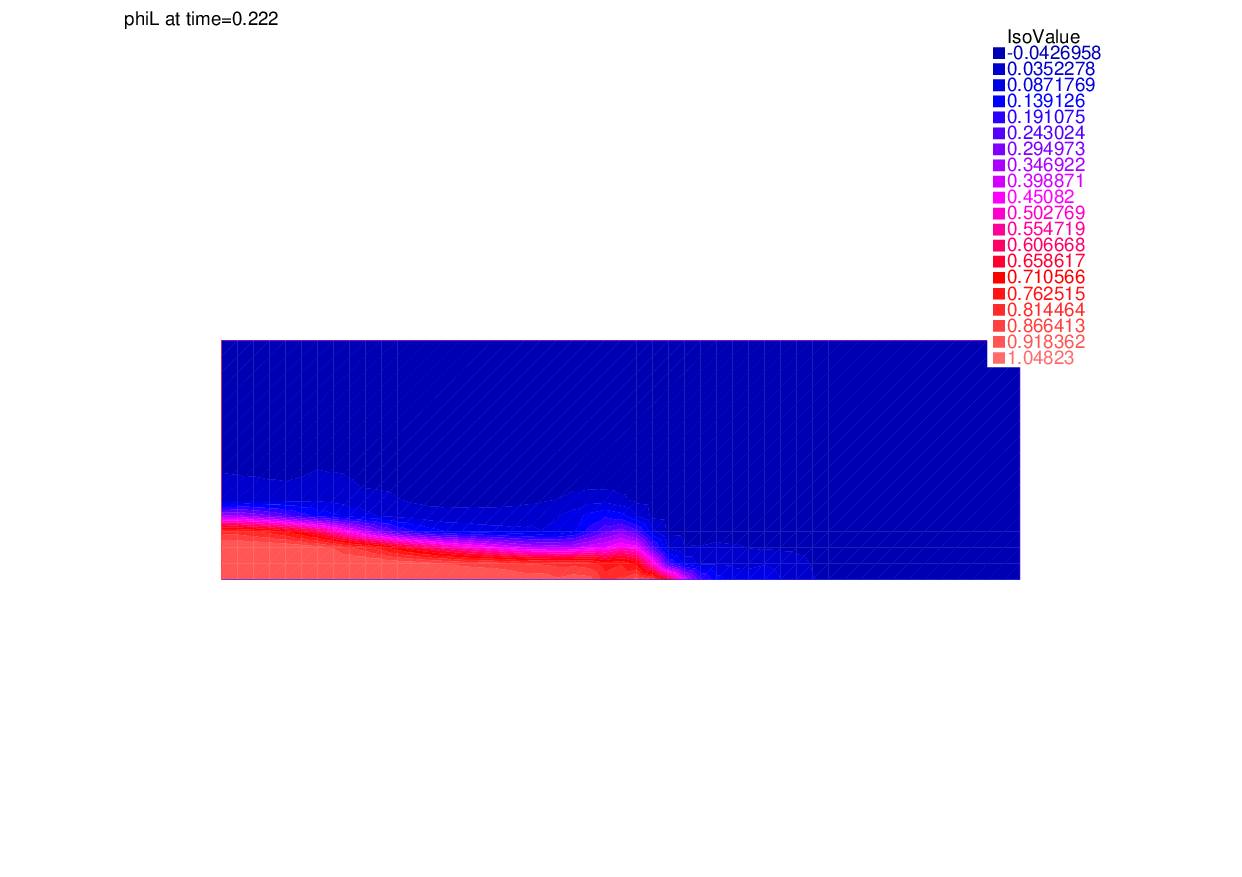}
		\caption{t=0.222(s)}
	\end{subfigure}
	\begin{subfigure}{0.49\textwidth}
		\centering
		\includegraphics[width=1.0\linewidth]{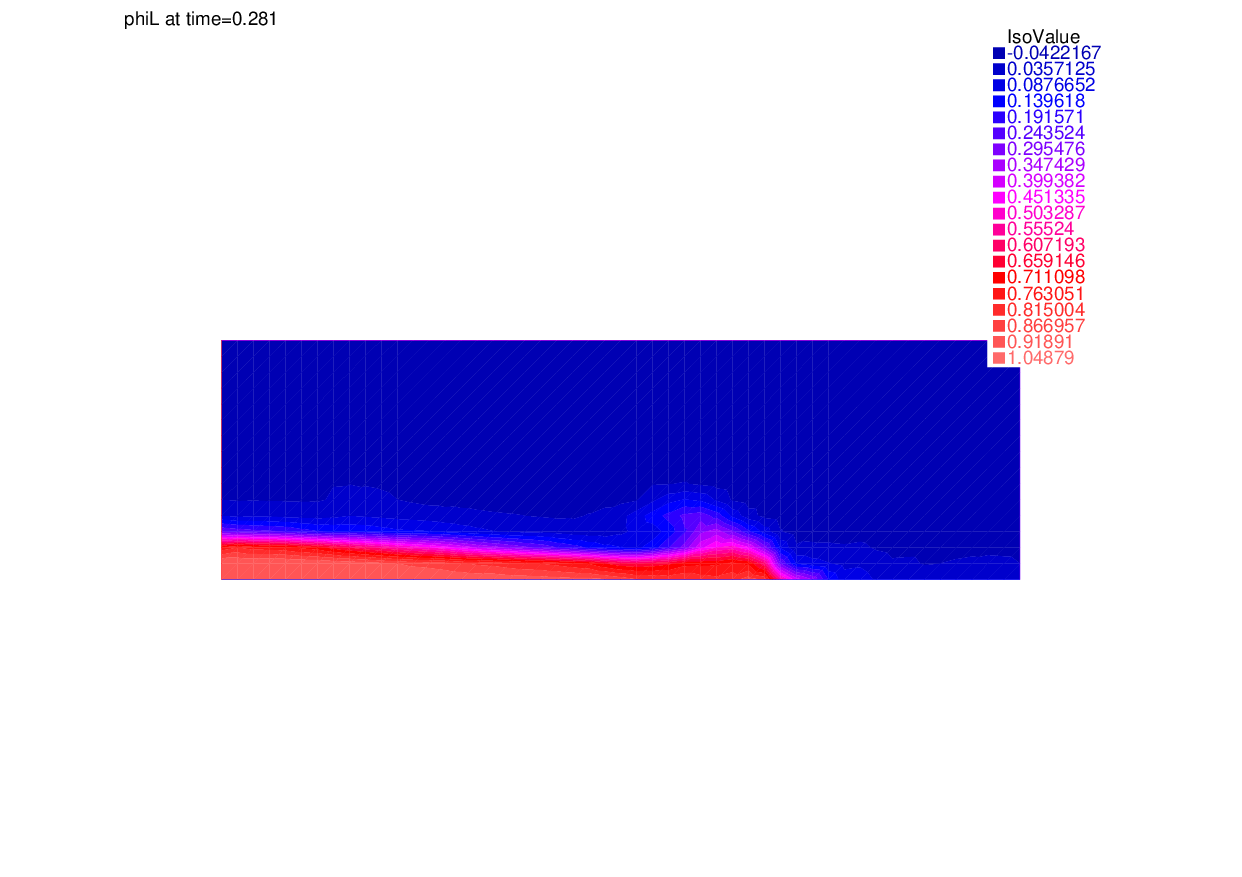}
		\caption{t=0.281(s)}
	\end{subfigure}
	
	\caption{Isovalues of the volume fraction for the dam break problem.}
	\label{fig:test2iso}
\end{figure}

\subsection{Bubble rising}
	This test simulates an air bubble rising under the effect of gravity in a closed box filled with water, as described in Figure \ref{test3domain}.  A $P_2$ element is used for approximating the velocities, while a $P_1$ element is used for all other variables. Similar to the case of the dam break problem, a hyperbolic tangent function is applied to smooth the interface between the bubble and the liquid. Furthermore, the setting of the drag force term is identical to that of the dam break problem to simulate the scenario where the bubble and its surrounding liquid are single-phase regions. The stabilization introduced in Eqs. \eqref{step5.1} and \eqref{step5.2} is adopted, just as in the dam break problem. Initially, both the air bubble and water are at rest, and the pressure follows a hydrostatic profile. For the computation, a non-uniform mesh depicted in Figure \ref{meshfig} is used to provide a high-quality representation of the initial conditions. The resulting profiles of the gas volume fraction at different physical times \(t = 0,~0.15,~0.35,~0.6,~0.75,~1\) are shown in Figure \ref{fig:test3iso}. The shape of the bubble after sufficient evolution, as depicted in Figure \ref{fig:test3iso}, resembles the case with a high Bond number presented in \cite{hua2007numerical}. However, at the final step, the interfaces fragmented into more pieces compared to the results in \cite{hua2007numerical}. This discrepancy might be attributed to the absence of a capillary effect model in our simulation, unlike the setup in \cite{hua2007numerical}, where the capillary model was applied (though not explicitly specified in that paper). Moreover, the use of stiffened models for equations of state in \cite{hua2007numerical} might lead to inaccuracies in the computed densities of the gas after it rises to a certain height, contributing to another source of discrepancy. Additionally, we observe that numerical diffusion does not inhibit the development of interface instabilities.

\begin{figure}
\centering
\includegraphics[width=0.5\textwidth]{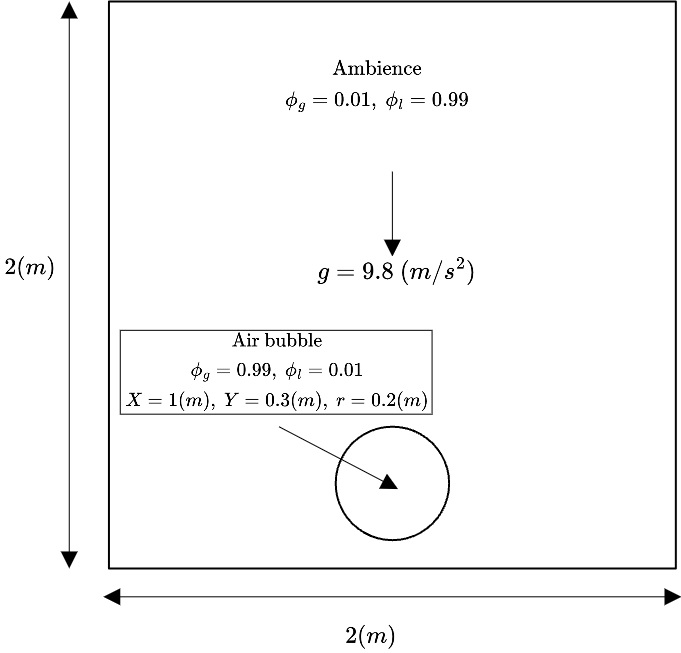}
\caption{Initial configuration of the bubble rising problem.}
\label{test3domain}
\end{figure}

\begin{figure}
\centering
\includegraphics[width=0.5\textwidth]{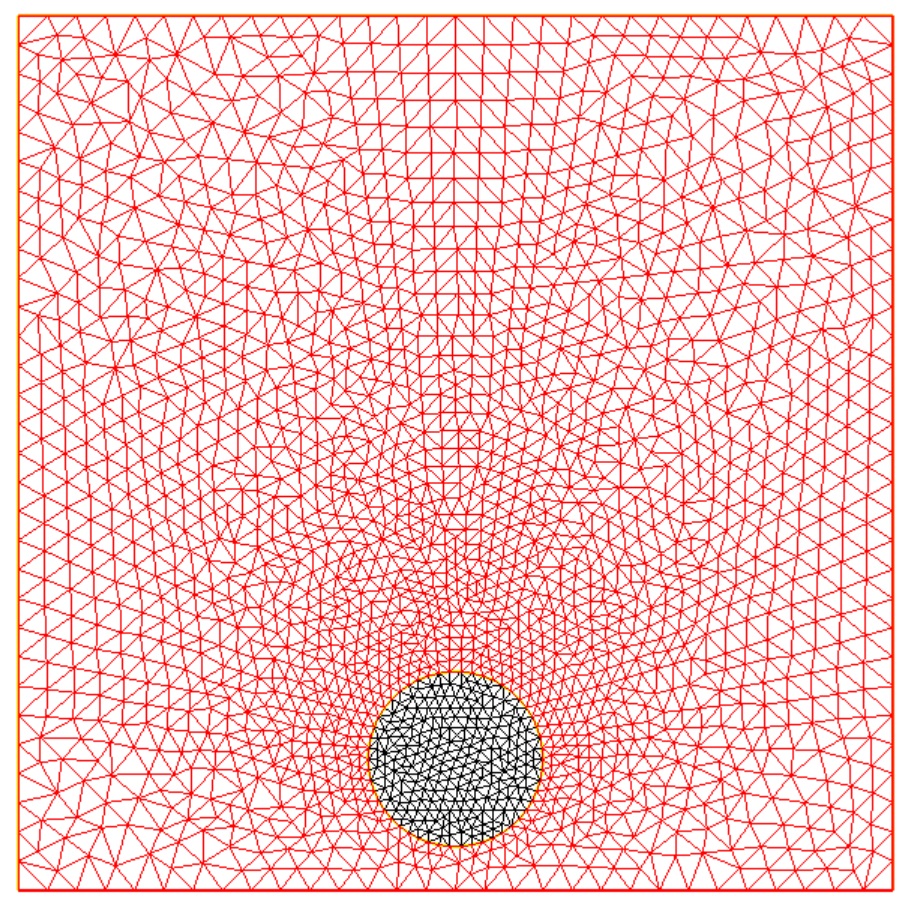}
\caption{The mesh in use for the simulation of bubble rising problem.}
\label{meshfig}
\end{figure}

\begin{figure}
\centering
	\begin{subfigure}{0.5\textwidth}
		\centering
		\includegraphics[width=1.0\linewidth]{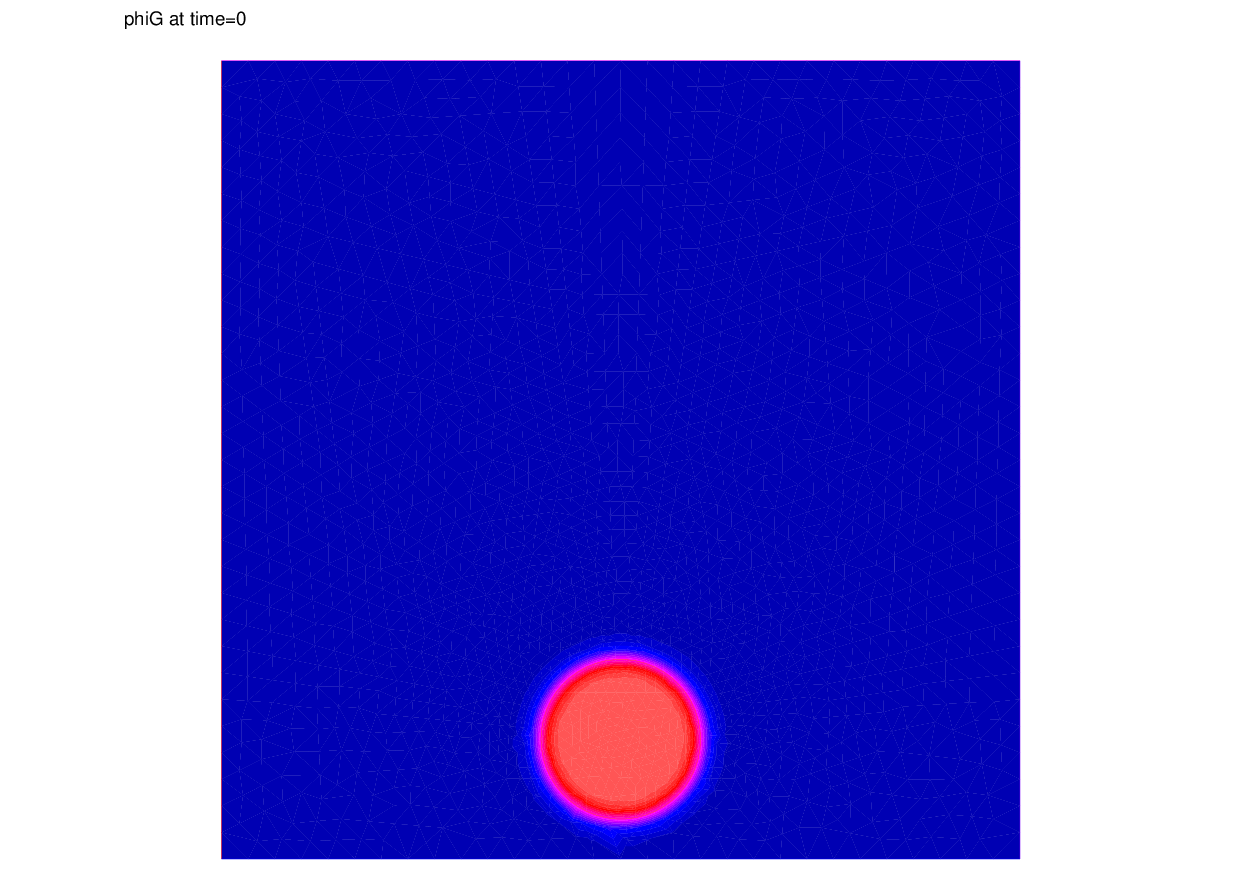}
		\caption{t=0(s)}
	\end{subfigure}\hspace{\fill}%
	\begin{subfigure}{0.5\textwidth}
		\centering
		\includegraphics[width=1.0\linewidth]{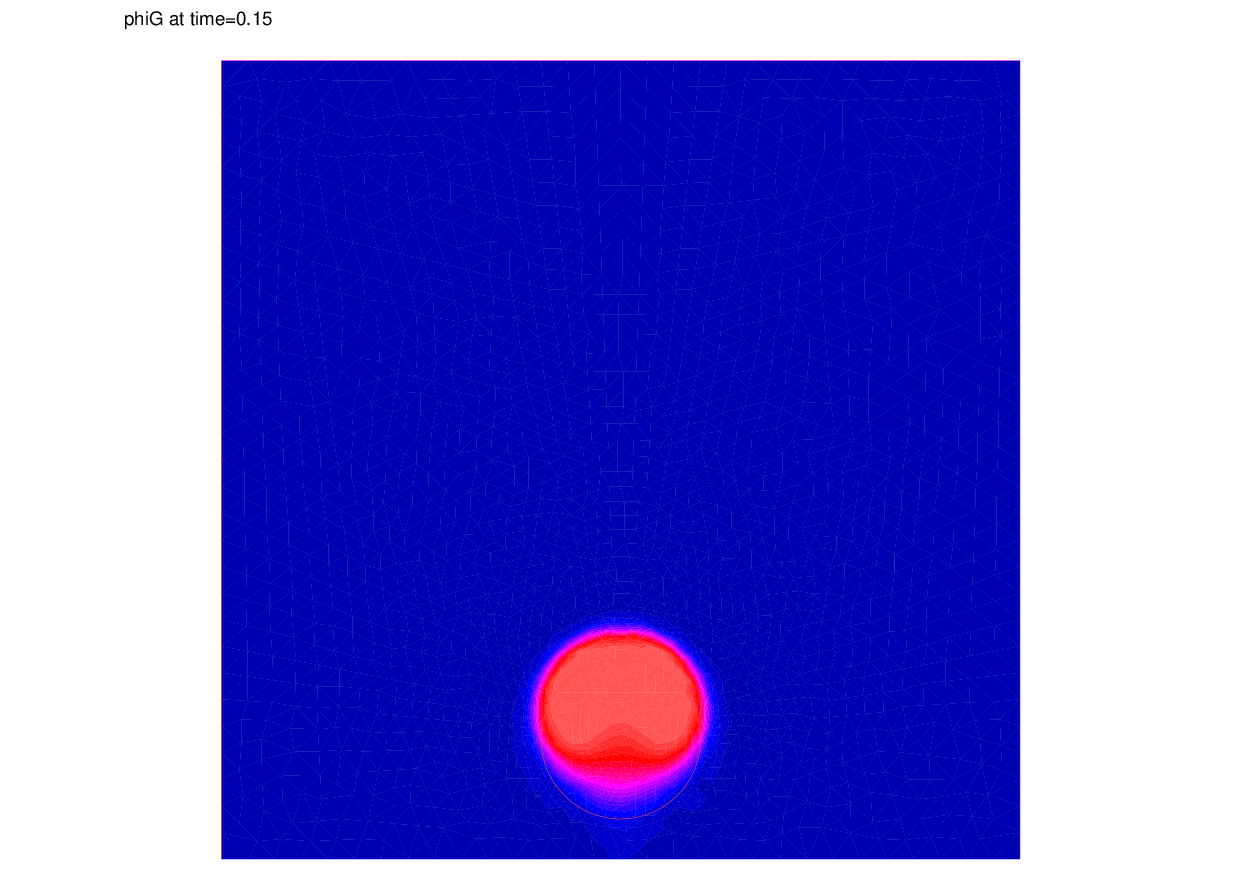}
		\caption{t=0.15(s)}
	\end{subfigure}
	\begin{subfigure}{0.5\textwidth}
		\centering
		\includegraphics[width=1.0\linewidth]{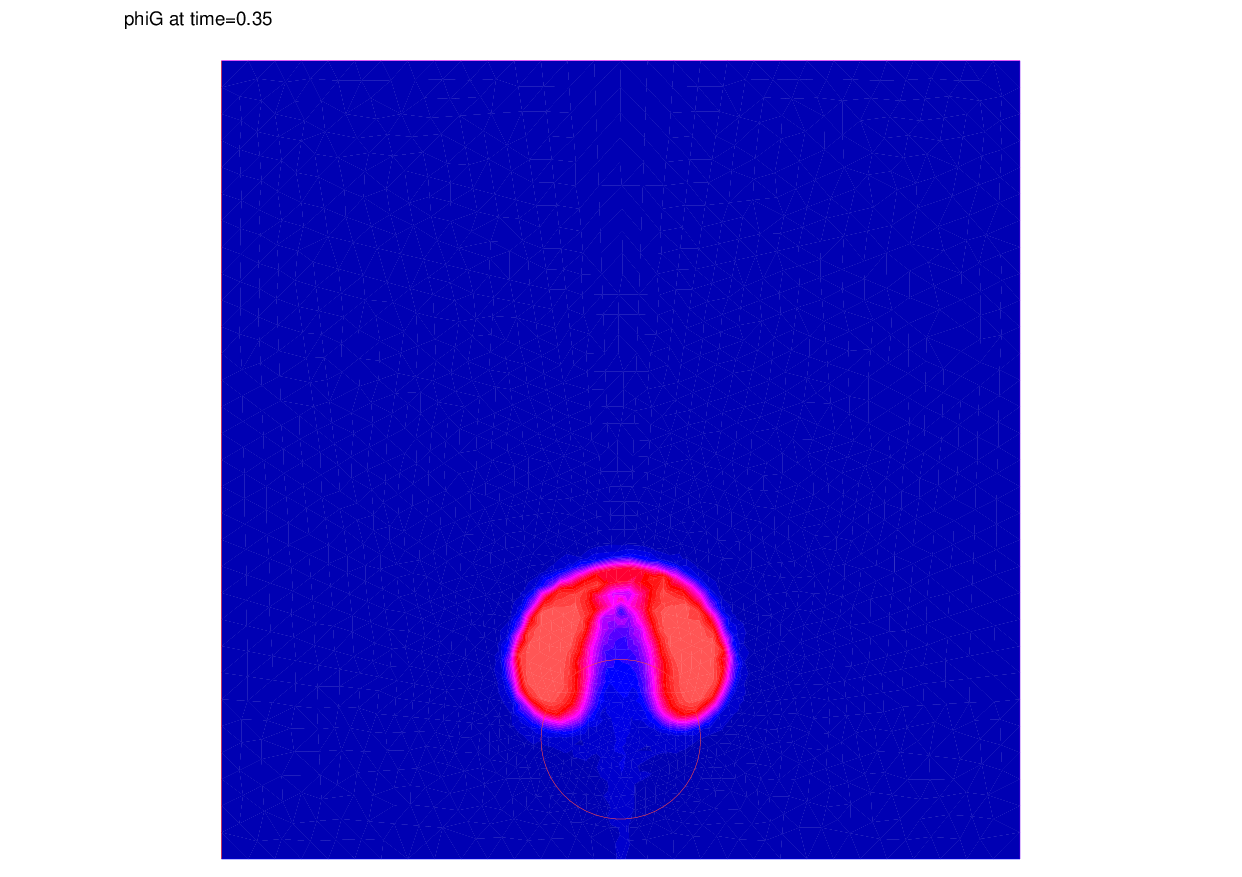}
		\caption{t=0.35(s)}
	\end{subfigure}\hspace{\fill}%
	\begin{subfigure}{0.5\textwidth}
		\centering
		\includegraphics[width=1.0\linewidth]{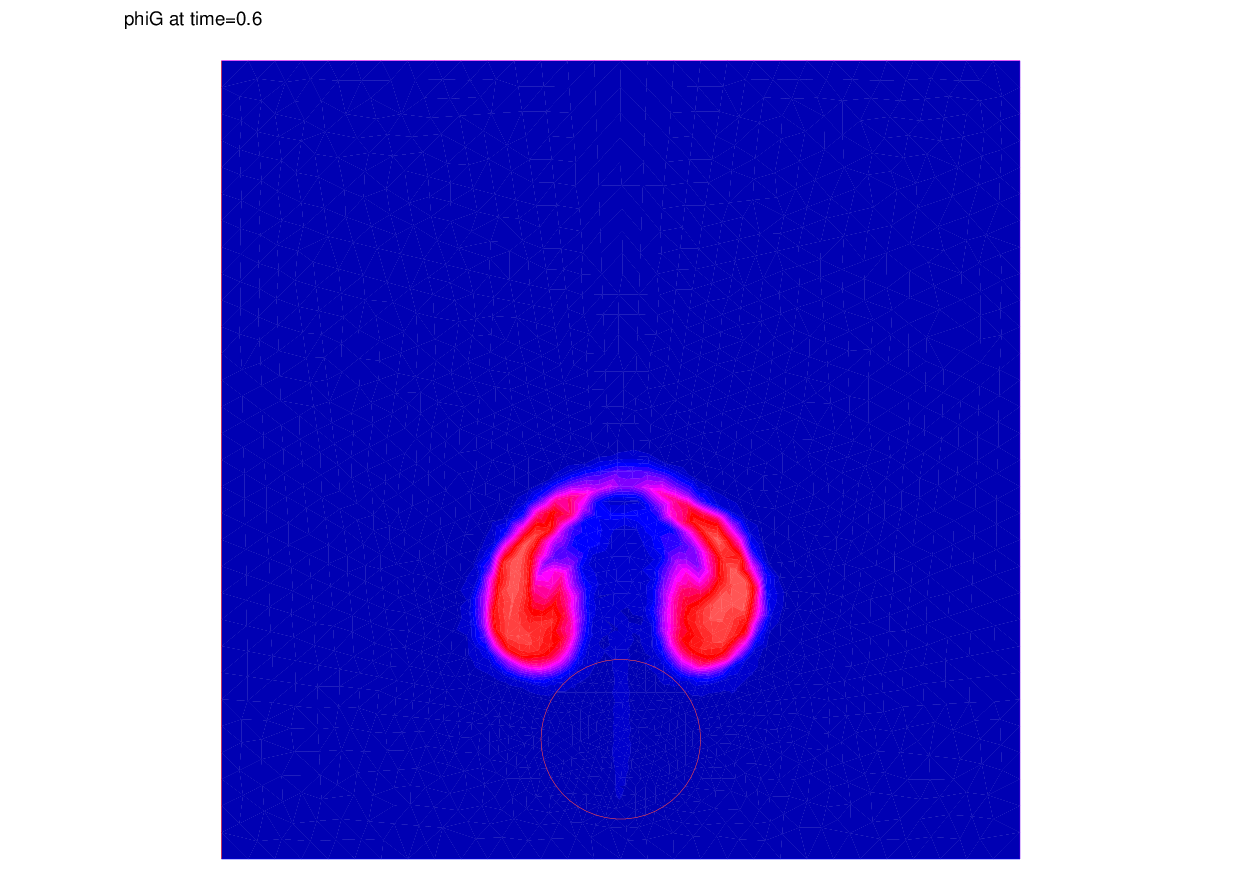}
		\caption{t=0.6(s)}
	\end{subfigure}
	\begin{subfigure}{0.49\textwidth}
		\centering
		\includegraphics[width=1.0\linewidth]{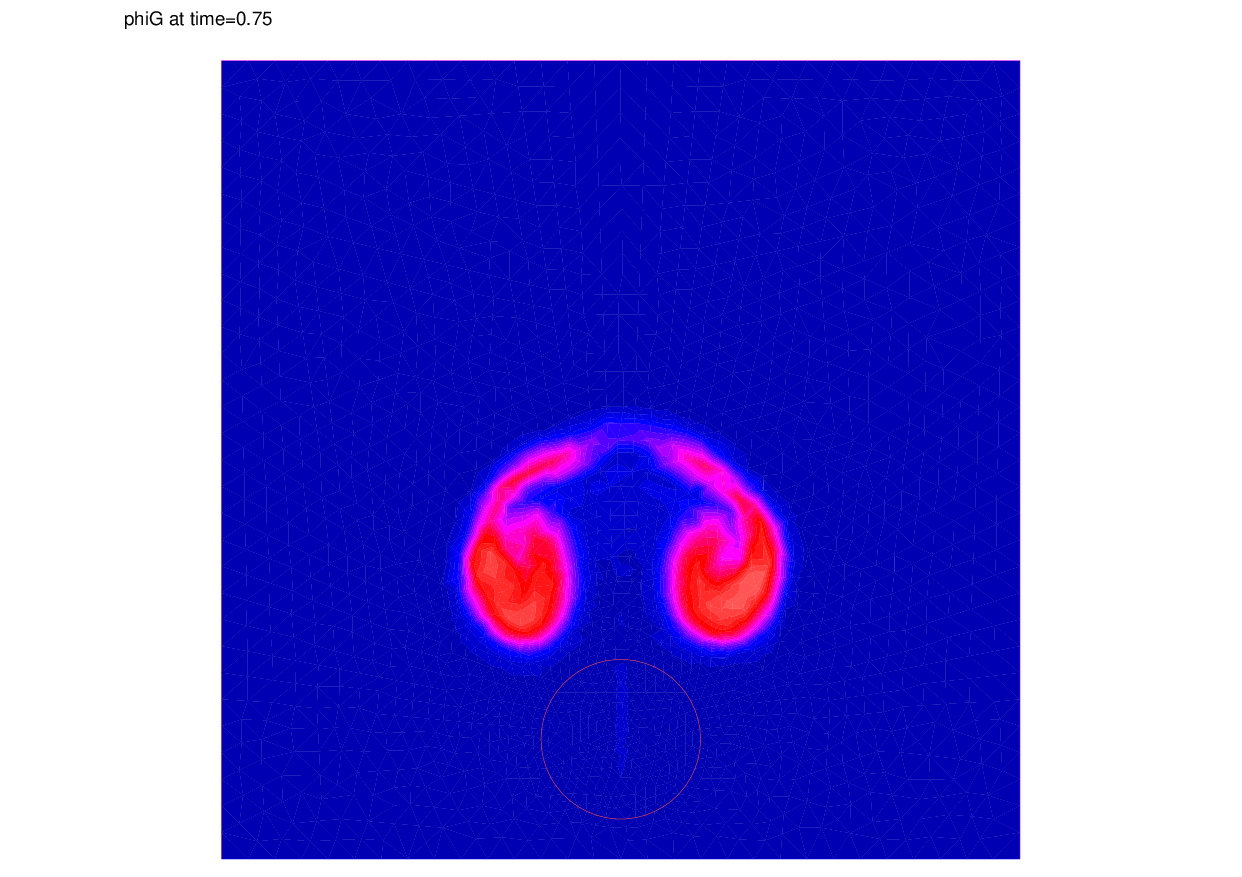}
		\caption{t=0.75(s)}
	\end{subfigure}
	\begin{subfigure}{0.5\textwidth}
		\centering
		\includegraphics[width=1.0\linewidth]{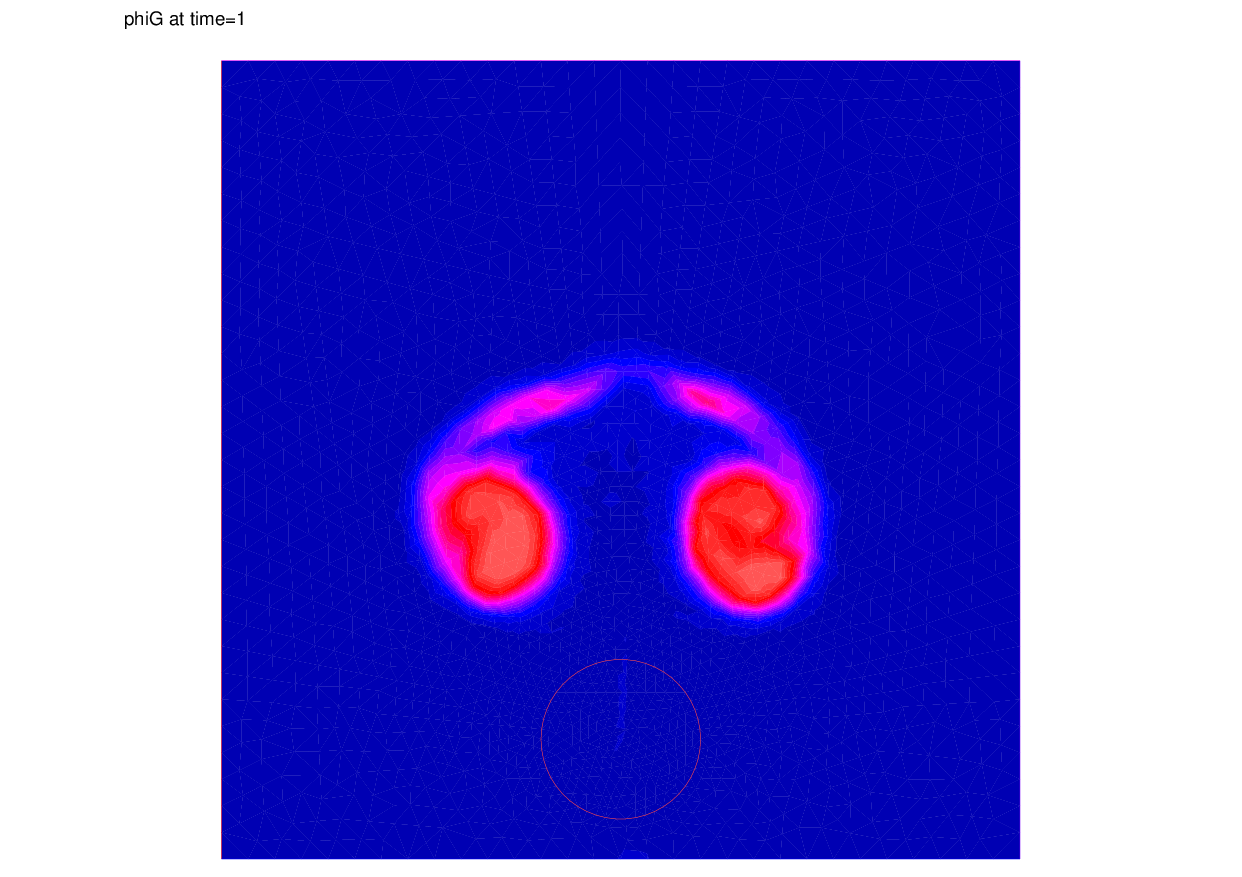}
		\caption{t=1(s)}
	\end{subfigure}
	\caption{Isovalues of the volume fraction for the bubble rising problem.}
	\label{fig:test3iso}
\end{figure}

\section{Conclusion}
In this work, a new projection method for solving a viscous two-fluid model is proposed. It is demonstrated that all iterations can be confined to the projection step, which holds promise for accelerating computations if an improved iteration technique for Step 5 is developed. Furthermore, this method is proven to be unconditionally stable at the time-discrete level. A suitable assignment of intermediate densities and volume fractions maintains the stability of the numerical scheme, as justified by the stability analysis for the time-discrete problem. A crucial stabilizing term is identified to address the instability caused by abnormal fluid divergence, which is validated by numerical experiments. The first-order temporal accuracy of the scheme is confirmed by the first numerical test in Section 6. Additionally, the proposed numerical scheme shows good agreement with some classical benchmarks.

\appendix

\bibliographystyle{siamplain}
\bibliography{ref}

\section*{Acknowledgement}
I would like to express my gratitude to Olivier Pironneau for his suggestions after reading this work, as well as to Li-Hsuan Shen and Chien-Ting Wu for their fruitful discussions. I also acknowledge the financial support from the National Science and Technology Council of Taiwan (NSTC 112-2639-E-011-001-ASP).

\section{Derivation of Eq. \eqref{goveq1}}
Most natural flows involve multiple fluids. This category includes the flow of immiscible fluids such as air and water, or oil, gas, and water, among other examples. When two fluids are miscible, they generally form a new fluid with distinct rheological properties. A stable emulsion of water and oil is a notable example: while both water and oil are Newtonian fluids, their emulsion behaves as a non-Newtonian fluid.

A simple yet non-trivial example of multi-fluid flow is the propagation of small-amplitude waves at the interface between a gas and a liquid, commonly referred to as separated flow in two-fluid terminology. In this scenario, each fluid obeys its own set of equations, with interactions occurring at the free surface. One mathematical model used to describe this interaction is the Navier-Stokes equation with a free surface: \begin{equation}\label{compns1} \partial_t \rho + \nabla\cdot (\rho \bm{u}) = 0, \end{equation} \begin{equation}\label{compns2} \partial_t(\rho\bm{u}) + \nabla\cdot(\rho\bm{u}\otimes\bm{u}) + \nabla p - \nabla\cdot \bm{\tau} = \rho \bm{g}, \end{equation} where $\bm{g}$ is the gravitational acceleration, $\rho$ is the fluid density, $\bm{u} = (u_x, u_y, u_z)^T$ is the fluid velocity vector, and $\bm \tau$ represents the stress tensor. The volume occupied by the two fluids (gas and liquid) is separated by a free surface $z=\eta(t,x,y)$, which satisfies: \begin{equation}\label{compns3} \partial_t \eta + u_x \partial_x \eta + u_y \partial_y \eta = u_z, \end{equation} where the region $z > \eta$ is occupied by the gas with density $\rho = \rho_g$, and the region $z < \eta$ is occupied by the liquid with density $\rho = \rho_l$. To close the system, it is assumed that each fluid follows its own equation of state: \begin{equation}\label{compns4} \mathcal{F}_k (\rho_k, p) = 0, \quad k=g,l. \end{equation}

As wave amplitude increases, wave breaking may occur. Near the gas-liquid interface, small liquid droplets can be suspended in the gas, while gas bubbles may form within the liquid. These inclusions can be very small. The collapse and fragmentation of the interface can lead to a topologically complex free surface, involving a wide range of length scales. In such cases, two-fluid models become not only relevant but essential. A powerful approach to model such phenomena is volume averaging, which leads to the development of a two-fluid model.
We denote by $\phi_g(x,t)$ the volumetric fraction of gas and by $\phi_l(x,t)$ that of liquid. This means that in an infinitesimal volume $dx^3$ around the point $x$, the volume occupied by gas is $\phi_g(x,t)dx^3$, while the volume occupied by liquid is $\phi_l(x,t)dx^3$. Therefore, we have $\phi_g + \phi_l = 1$.
Applying the volume averaging procedure (see, e.g., \cite{Ni1991, ishii2010thermo}) to \eqref{compns1} and \eqref{compns2}, and neglecting some fluctuation terms, we obtain the following equations:
\begin{subequations}\label{goveq1e}
\begin{equation}
\phi_g + \phi_l = 1, \quad \text{in}~\Omega_T,
\end{equation}    
\begin{equation}
\partial_t(\phi_k \rho_k) + \nabla\cdot(\phi_k \rho_k \bm{u}_k) = 0, \quad \text{in}~\Omega_T,
\end{equation}
\begin{equation}\label{goveq1.3e}
\partial_t (\phi_k \rho_k \bm{u}_k) + \nabla\cdot(\phi_k \rho_k \bm{u}_k \otimes \bm{u}_k) - \nabla\cdot(\phi_k \tau_k(\bm{u}_k)) 
+ \phi_k \nabla p_k + F_{D,k} = \phi_k \rho_k \bm{g}, \quad \text{in}~\Omega_T,
\end{equation}
\end{subequations}
where $F_{D,k}$ are drag forces satisfying certain constitutive laws and $\tau_k$ are viscous stresses. \vfill
\end{document}

%% file: tmp_docsiamart_header.tex
\title{A Projection Method for Compressible Generic Two-Fluid Model\thanks{Submitted to the editors DATE.
\funding{I acknowledge the financial support from the National Science and Technology Council of Taiwan (NSTC 112-2639-E-011-001-ASP).}}}

\author{Po-Yi Wu\thanks{Sustainable Electrochemical Energy Development (SEED) Center, National Taiwan University of Science and Technology, Taipei City 106, Taiwan (\email{r04527030@gmail.com}).}}

\headers{A Projection Method for Compressible Generic Two-Fluid Model}{Po-Yi Wu}

%% file: tmp_docsiamart_abstract.tex
\begin{abstract}
A new projection method for a generic two-fluid model is presented in this work. Specifically, we extend the projection method, originally designed for single-phase variable density incompressible and compressible flows, to viscous compressible two-fluid flows. The key idea is that the single pressure $p$ can be uniquely determined by the products of volume fractions and densities $\phi_k \rho_k$, of the two fluids. Additionally, the stability of the method is ensured by appropriately assigning intermediate step variables at the time-discrete level and incorporating a stabilizing term at the fully discrete level. We prove the energy stability of the proposed numerical scheme, and its validity is demonstrated through three numerical tests.
\end{abstract}

\begin{keywords}
Two-fluid model, Projection method, Stability analysis
\end{keywords}

\begin{MSCcodes}
65M12, 65M60
\end{MSCcodes}